\documentclass[twocolumn,twoside]{IEEEtran}
\usepackage{amsmath,amssymb,euscript ,yfonts,psfrag,latexsym,dsfont,graphicx,bbm,color,amstext,wasysym,subfig,flushend,parskip}
\usepackage{times,mathptmx,epsfig,graphicx,color}
\usepackage{amsthm}
\usepackage{soul}

\graphicspath{{./},{./figures/},{../figures/},{../}}
\usepackage{graphicx}
\usepackage{url}
\usepackage{amsmath}
\usepackage{amsthm}
\usepackage{amsfonts}
\usepackage{array}
\usepackage{cancel}
\usepackage{enumitem}
\newcommand{\E}{{\mathbb E}}
\newcommand{\R}{{\mathbb R}}
\newtheorem{propo}{\textbf{Proposition}}

\newtheorem{probl}{\textbf{Problem}}
\newtheorem{rem}{\textbf{Remark}}
\newtheorem{defi}{\textbf{Definition}}
\newtheorem{thm}{\textbf{Theorem}}
\theoremstyle{definition}

\newtheorem{assume}{\textbf{Assumption}}
\usepackage{soul}
\usepackage{tabularx}

\newcommand{\black}{\color{black}}
\newcommand{\ignore}[1]{}
\DeclareMathOperator*{\trace}{trace}

\newcommand{\mathcalO}{O}
\newcommand{\bhs}{\hspace*{-2pt}}

\def\spacingset#1{\def\baselinestretch{#1}\small\normalsize}
\spacingset{1.02}

\begin{document}
\title{{
Regularized transport
%Entropic interpolation 
between\\
singular covariance matrices
%degenerate Gaussian distributions
\thanks{Partial support was provided by NSF under grants %1665031,
1807664, 1839441 (TTG) and 1901599 (YC), and by AFOSR under grant FA9550-17-1-0435 (TTG), and by the University of Padova Research Project CPDA 140897 (MP).}}
}
\author{Valentina Ciccone$^\dagger$, Yongxin Chen$^\ast$, Tryphon T.\ Georgiou$^\ddagger$, Michele Pavon$^\star$%
\thanks{$^\dagger$ Department of Information Engineering, University of Padova, Padova 35121, Italy,
{\small valentina.ciccone@dei.unipd.it}}%
\thanks{$^\ast$ School of Aerospace Engineering, Georgia Institute of Technology, Atlanta, GA 30332,
{\small yongchen@gatech.edu}}%
\thanks{$^\ddagger$ Department of Mechanical and Aerospace Engineering,
University of California, Irvine, CA 92697, USA,
{\small tryphon@uci.edu}}%
\thanks{$^\star$ Dipartimento di Matematica ``Tullio Levi-Civita",
Universit\`a di Padova, 35121 Padova, Italy,
{\small pavon@math.unipd.it}}%
}

\maketitle
\begin{abstract} We consider the problem of  %minimum-energy 
steering 
%of 
a linear stochastic system between two end-point \textit{degenerate} %multivariate 
Gaussian distributions in finite time. This accounts for those situations in which some but not all of the state entries are uncertain at the initial, $t=0$, and final time, $t=T$.
This problem entails  non-trivial technical challenges as the singularity of terminal state-covariance causes the control to grow unbounded at the final time $T$.  Consequently, the entropic interpolation ({\em Schr\"odinger Bridge}) is provided by a diffusion process which is not {\em finite-energy}, thereby placing this case outside of most of the current theory.
In this paper, we show that a feasible interpolation can be derived as a limiting case of earlier results for non-degenerate cases, and that it can be expressed in closed form. Moreover,
we show that such interpolation belongs to the same \textit{reciprocal
class} of the uncontrolled evolution. By doing so we also highlight a time-symmetry of the problem, contrasting dual formulations in the forward and reverse time-directions, where in each the control grows unbounded as time approaches the end-point (in the forward and reverse time-direction, respectively).
\end{abstract}


\section{Introduction} 
The problem of optimally steering a Markov process between two end-point marginal distributions has its roots in the  thought experiment of {\em large deviations} for independent Brownian particles formulated by Schr{\"o}dinger in the early 30's 
\cite{schrodinger1931} \cite{schrodinger1932theorie}. 
Since then, considerable literature has flourished and important connections have been  made bridging the Schr{\"o}dinger Problem with Stochastic Control Theory, \cite{dai1991stochastic, dai1990markov,chen2016relation} and, more recently, with Optimal Mass Transport, \cite{Mik,MT2006,MT2008,leonard2012schrodinger, leonard2013survey,gentil2015analogy}.

A special instance of this problem, which is of great interest from a control engineering perspective, is the one of steering a linear stochastic system between two end-points Gaussian distributions in finite time and with minimum energy \cite{beghi1997continuous, chen2016optimal2,chen2017optimal,
chen2018optimal,bakolas2016optimal,
halder2016finite}. In particular, such a problem can be recast as a {\em finite-horizon} covariance control problem in the spirit of the seminal works \cite{hotz1987covariance,collins1987theory}.
In \cite{chen2016optimal}, this problem has been studied for the general  case of possibly degenerate diffusions and the solution for the optimal control problem has been   explicitly
derived  by  solving  two  differential  Lyapunov   equations  non linearly coupled through their boundary conditions.

The purpose of the current paper is to study in detail the case in which the desired marginals are Gaussian with singular covariances and, consequently, have no density. %In fact, this apparently peripheral instance can be seen to be of great interest from both an applied and a theoretical perspective.
This leads to non trivial technical issues  
as 
%the final distribution being {\vale not absolutely continuous} with respect to the reference measure causes 
the control becomes unbounded at the terminal time $T$. It follows that the bridge process solving the problem is no more a {\em finite-energy diffusion}. But this is the key assumption under which  most of the theory has been developed \cite{F2}. 
Our goal is to show that a feasible interpolation can be derived as a limiting case of the results in \cite{chen2016optimal} and that it can be expressed in closed form. Moreover, we aim to show that this interpolation belongs to the same {\em reciprocal class} of the uncontrolled evolutions, 
namely, the class of probability laws having the same three-point transition densities, see \cite{Jam1,LK}; this is the case for nondegenerate marginals \cite{F2}.
As a by-product of our analysis, we also gain interesting insight into the time-symmetry of the problem.

The paper is organized as follows. In Section \ref{BrownianB}, as exemplification for our studies, we review some basic facts about the Brownian Bridge. In Section \ref{Review} we recall some results from \cite{chen2016optimal}. Then, in Section \ref{Problem_statement}, we formally state our problem whereas in Section \ref{Results} we present our main results. In Section %\ref{time-symmetry}, we briefly comment on the time-symmetry of the solution. Finally, in Section 
\ref{Numerical_ex} a numerical example is worked out for illustrative purposes. Finally, some open points and future directions are discussed in Section \ref{Conclusion}. The less relevant proofs are deferred to the Appendix.

%We are interested in the problem of steering a particle systems between two end-point marginal distributions in finite time and with minimum energy. Such problem is closely related to the so called Schroedinger Bridge Problem (SB) and to Optimal Mass Transport (OMT). Both of these problem admits a stochastic control formulation.

\section{The Brownian Bridge}\label{BrownianB}
To motivate our successive analysis and illustrate some of the difficulties we are going to face, we briefly recall some basic facts about a well-known example of infinite-energy diffusion process, namely the Brownian Bridge.
The {\em Brownian Bridge} is  defined as the continuous-time stochastic process $(X_t)_{t\in[0,1]}$ which is obtained from a standard $n$-dimensional Wiener process $(W_t)_{t\in[0,1]}$ by conditioning on $W_1=0$, \cite{RY1999}, \cite{chen2016stochastic}. Hence both the initial and final marginals are Dirac delta masses concentrated at $0\in\R^n$. The Brownian bridge satisfies the following  stochastic differential equation
\begin{align}\label{BB_sde}
dX_t=-\frac{X_t}{1-t}dt+dW_t, \quad 0\leq t < 1, \; X_0=0.
\end{align}
%Then, as an academic exercise, we want to show that the process $X_t$, which is the strong solution of \eqref{BB_sde} in $[0,1)$, is indeed the sought bridge.By standard computation 
A straightforward calculation gives
$$X_t=(1-t)\int_0^{t}\frac{1}{1-s}dW_s.$$
$X_t$ is clearly a zero-mean Gaussian process. 
%Thus, it only remains to show that $$L^{2}-\lim_{t \nearrow 1} X_t=0.$$
Let $P(t):=\mathbb{E}[X_tX_t']$ denote the variance of $X_t$. By a standard computation, $P(t)=t(1-t)I_n$. Then, since $\lim_{t\nearrow 1} P(t)=0$, we conclude that $X_t$ converges to zero in mean square. 
%$L^{2}(\Omega,\mathcal{F},\lbrace \mathcal{F}_t\rbrace,\mathbb{P})$. 
We compute now
\begin{align*}
&\trace\E\left\{\int_0^1\left(\frac{-X_t}{1-t}\right)\left(\frac{-X_t}{1-t}\right)'dt\right\}\\
&=\trace\int_0^1P(t)\frac{1}{(1-t)^2}dt
=n\int_0^1\frac{t}{1-t}dt=+\infty.
\end{align*}
Thus, the Brownian bridge is {\em not} a finite-energy diffusion.

\section{Review of Previous Results}\label{Review}
Let $(\xi_t)_{t\in [0,T]}$ be an $\mathbb{R}^{n}$-valued stochastic diffusion process defined on $(\Omega, \mathcal{F},\lbrace \mathcal{F}_t \rbrace, \mathbb{P})$ and satisfying the linear stochastic differential equation
\begin{align}
\label{sde_pI}
%\label{prior}
d\xi_t=A(t)\xi_t dt+B(t)dW_t, \quad \xi_{t=0}=\xi_0 \; \text{a.s.}
\end{align}
where $(W_t)_{t\in [0,T]}$ is a standard m-dimensional $\mathcal{F}_t$-Wiener process, $A(t):[0,T]\rightarrow \mathbb{R}^{n\times n}$ and $B(t):[0,T]\rightarrow \mathbb{R}^{n\times m}$ are bounded and continuous matrix functions and $\xi_0$ is an n-dimensional random vector independent of $(W_t)_{t\in [0,T]}$, $\xi_0\sim \rho_0$ with $\rho_0$ a zero-mean Gaussian distribution with covariance $\tilde{\Sigma}_0$.

The problem of forcing the diffusion process $(\xi_t)_{t\in[0,T]}$ to a desired end-point probability distribution $\rho_T$, with $\rho_T$ zero-mean Gaussian with covariance $\tilde{\Sigma}_T$,
has been considered in \cite{chen2016optimal}. To this end, the controlled diffusion process  
$(\xi_t^{u})_{t\in [0,T]}$
\begin{align}
\label{controlled_sde_pI}
d\xi^{u}_t=A(t)\xi^{u}_tdt+B(t)u(t)dt+B(t)dW_t, \quad \xi^{u}_{t=0}=\xi_0 \; \text{a.s.}
\end{align}
is introduced, where $u\in\tilde{\mathcal{U}}$, the class of admissible controls specified as follows.
\begin{defi}
A control $u(\cdot)\in\tilde{\mathcal{U}}$ if
 $u(t)$ is $\mathcal{F}_t$-adapted and if, in addition,
the stochastic differential equation \eqref{controlled_sde_pI} admits a strong solution in $[0,T]$ and
 $\mathbb{E}\big\{\int_0^{T}u(t)'u(t)dt\big\}<\infty$.
\end{defi}
Our problem can be formally stated as follows:
\begin{probl}\label{problem_pI} Provided that the set $\tilde{\mathcal{U}}$ is non-empty, find an admissible control $u^{*}(\cdot)\in\tilde{\mathcal{U}}$ such that:
1) 
$\xi^{u^{*}}_{t=0}$ is distributed according to $\rho_0$, and $\xi^{u^{*}}_{t=T}$ according to $\rho_T$, and
2) among all the admissible controls satisfying the previous point, $u^{*}$ minimizes the cost function 
$$ J(u):=\mathbb{E}\Big\{\int_0^{T}u(t)'u(t)dt\Big\}.$$ 
\end{probl}
A solution was given in \cite{chen2016optimal} under the following assumptions.
\begin{assume}\label{a1}
The covariances $\tilde{\Sigma}_0,\;\tilde{\Sigma}_T$ are positive definite.
\end{assume}
\begin{assume}\label{a2}
The pair $(A(t),B(t))$ is controllable, in the sense that the reachability Gramian
$$ M(t_1, t_0):=\int_{t_0}^{t_1}\Phi(t_1, \tau)B(\tau)B(\tau)'\Phi(t_1,\tau)'d\tau$$
is non singular for all $t_0<t_1$, $t_0,t_1\in [0,T]$, where $\Phi(t,s)$ denotes the system state-transition matrix of $A(\cdot)$.
%determined via
%$$\dfrac{\partial}{\partial t}\Phi(t,s)=A(t)\Phi(t,s), \quad \text{with }\; \Phi(t,t)=I$$
%which is non-singular for all $t,s \in [0,T]$.
\end{assume}

The next two results, Propositions 1 and Theorem 1 (from \cite{chen2016optimal}), are the starting point of our analysis.

\begin{propo}[\cite{chen2016optimal}]\label{propo_pI} Under Assumptions \ref{a1} and \ref{a2}, the following system of differential Lyapunov equations 
\begin{subequations}\label{system_Lyap}
\begin{align}\label{system_Lyap_P}
\dot{P}(t)&=A(t)P(t)+P(t)A(t)'+B(t)B(t)'\\ \label{system_Lyap_Q}
\dot{Q}(t)&=A(t)Q(t)+Q(t)A(t)'-B(t)B(t)'
\end{align}
\end{subequations}
coupled through the boundary conditions
\begin{subequations}\label{boundary_Lyap}
\begin{align}
\tilde{\Sigma}_0^{-1}&=P(0)^{-1}+ Q(0)^{-1}\\
\tilde{\Sigma}_T^{-1}&=P(T)^{-1}+Q(T)^{-1}.
\end{align}
\end{subequations}
admits a unique solution pair $(P(\cdot), Q(\cdot))$  such that $Q(t)$ and $P(t)$ are both non-singular on $[0,T]$. The pair is determined by \eqref{system_Lyap} and the initial conditions\footnote{We use the short hand notation $\Phi=\Phi(T,0)$, $M=M(T,0)$.}
\begin{small}
\begin{align*}
Q(0)=& \tilde{\Sigma}_0^{1/2} \Bigg( \dfrac{1}{2}I+\tilde{\Sigma}_0^{1/2}\Phi'M^{-1}\Phi \tilde{\Sigma}_0^{1/2} \\ 
& \hspace*{-0.0cm}-\left( \dfrac{1}{4}I+\tilde{\Sigma}_0^{1/2}\Phi'M^{-1}\tilde{\Sigma}_T M^{-1}\Phi \tilde{\Sigma}_0^{1/2} \right)^{1/2} \Bigg)^{-1} \tilde{\Sigma}_0^{1/2}\\
P(0)=& (\tilde{\Sigma}_0^{-1}-Q(0)^{-1})^{-1}.
\end{align*}\end{small}
\end{propo}

\begin{thm}[\cite{chen2016optimal}]\label{Thm_part_I} Under Assumptions \ref{a1} and \ref{a2}, Problem \ref{problem_pI} admits a unique optimal solution
$$u^{*}(t)=-B(t)'Q(t)^{-1}\xi_t$$
where $Q(t)$, together with the corresponding $P(t)$, solves the system of Lyapunov differential equations in Proposition \ref{propo_pI}.
\end{thm}

We state a dual to Proposition \ref{propo_pI} that will be used later on; the proof of Proposition \ref{rem_on_P(T)} is deferred to the Appendix.\\

\begin{propo}\label{rem_on_P(T)}
The solution $(P(\cdot),Q(\cdot))$ in Proposition \ref{propo_pI} can be equivalently specified by system \eqref{system_Lyap} and the final conditions:
\begin{small}
\begin{align*}
P(T)=& \tilde{\Sigma}_T^{1/2} \Bigg( \dfrac{1}{2}I+\tilde{\Sigma}_T^{1/2}M^{-1}\tilde{\Sigma}_T^{1/2} \\ 
& \hspace*{-0.0cm} -\left( \dfrac{1}{4}I+\tilde{\Sigma}_T^{1/2}M^{-1}\Phi\tilde{\Sigma}_0 \Phi'M^{-1} \tilde{\Sigma}_T^{1/2} \right)^{1/2} \Bigg)^{-1} \tilde{\Sigma}_T^{1/2}\\
Q(T)=& (\tilde{\Sigma}_T^{-1}-P(T)^{-1})^{-1}.
\end{align*}
\end{small}
%\blue{\st{where again we have used the short hand notation $\Phi=\Phi(T,0)$, $M=M(T,0)$.}}
\end{propo}

\section{Problem statement}
\label{Problem_statement}
Let $(\zeta_t)_{t\in [0,T]}$ be a $\mathbb{R}^{n}$-valued stochastic diffusion process defined on $(\Omega, \mathcal{F},\lbrace \mathcal{F}_t \rbrace, \mathbb{P})$ and 
satisfying the linear stochastic differential equation
\begin{align}
\label{system}
d\zeta_t=A(t)\zeta_tdt+B(t)dW_t, \quad \zeta_{t=0}=\zeta_0 \; \text{a.s.}
\end{align}
where %$(W_t)_{t\in [0,T]}$ is a standard m-dimensional Wiener process defined on $(\Omega, \mathcal{F}, \mathbb{P})$, $A(t):[0,T]\rightarrow \mathbb{R}^{n\times n}$ and $B(t):[0,T]\rightarrow \mathbb{R}^{n\times m}$ are (bounded and)  continuous matrix functions, while
$(W_t)_{t\in [0,T]}$, $A(\cdot)$ and $B(\cdot)$ are as above and $\zeta_0$ is an $n$-dimensional, zero-mean random vector independent of $(W_t)_{t\in [0,T]}$. We suppose that $\zeta_0$ is distributed according to a \textit{degenerate} multivariate normal $\mu_0$ with zero mean and singular covariance $\Sigma_0$, with $\text{rank}(\Sigma_0)=k_0<n$. %Then $\zeta_0$ has no density  w.r.t. the Lebesgue measure in $\mathbb{R}^{n}$,
% whose density $\mu_0$ is singular w.r.t. the Lebesgue measure in $\mathbb{R}^{n}$.
%however, by disintegration we can define a restriction of the Lebesgue measure to the $k_0$-dimensional affine subspace where the distribution is supported. Then, w.r.t. this measure, $\zeta_0$ has density
%\[\mu_0(\zeta)= (\text{det}^{*}(2\pi\Sigma_0))^{-1/2}\exp \left(-\dfrac{1}{2}\zeta'\Sigma_0^{\dagger}\zeta\right) \]
%where ${\dagger}$ denotes the generalized inverse and $\det^{*}$ is the pseudo determinant.
We consider the problem of steering system \eqref{system} to a desired final distribution $\mu_T$. We assume that $\mu_T$ is also a \textit{degenerate} multivariate normal with zero mean and singular covariance matrix $\Sigma_T$, with $\text{rank}(\Sigma_T)=k_T<n$.
Clearly, both $\mu_0$ and $\mu_T$ are not absolutely continuous with respect to the Lebesgue measure in $\mathbb{R}^{n}$ and therefore have no density with respect to such measure.
%which is also assumed to be {\vale not absolutely continuous} with respect to the Lebesgue measure in $\mathbb{R}^{n}$. %As before, by disintegration of measure we consider the restriction of the Lebesgue measure to the $k_T$-dimensional affine subspace of $\mathbb{R}^{n}$ where $\mu_T$ is supported and we define $\mu_T$ we respect to such measure as
%\[\mu_T(\zeta)=(\text{det}^{*}(2\pi\Sigma_T))^{-1/2}\exp \left(-\dfrac{1}{2}\zeta'\Sigma_T^{\dagger}\zeta\right) \]
%where $\Sigma_T\in\mathbb{R}^{n\times n}, \Sigma_T\succeq 0 , \; \text{rank}(\Sigma_T)=k_T$.\\
We now seek to formulate an atypical stochastic control problem with a two-point boundary conditions.
The corresponding controlled evolution is $(\zeta_t^{u})_{t\in [0,T]}$ satisfying
\begin{align}
\label{controlled_ev}
d\zeta^{u}_t=A(t)\zeta^{u}_tdt+B(t)u(t)dt+B(t)dW_t, \,\, \zeta^{u}_{t=0}=\zeta_0 \; \text{a.s.}
\end{align}
where $u\in\mathcal{U}$, the family of admissible control for the problem at hand, defined below.

\begin{defi} A control $u(\cdot)\in\mathcal{U}$ if
$u(t)$ is $\mathcal{F}_t$-adapted and if, in addition, for any $\varepsilon\in(0,T)$, 
equation \eqref{controlled_ev} admits a strong solution in $[0,T-\varepsilon]$ and
 $u$ has finite energy on $[0,T-\varepsilon]$.
\end{defi}

We stress that \textit{any} admissible control steering the system to the desired %terminal 
degenerate marginal $\mu_T$, necessarily needs
to have infinity energy,
see\footnote{
 A Markov finite energy diffusion has both forward $b_+$ and reverse-time $b_-$ drifts with finite energy \cite{F2}.  The two are related via the density $\rho$ as in  $b_-=b_+-BB'\nabla\log\rho$. Thus, if $\rho$ tends to become singular for $t\nearrow T$, the energy of at least one of the drifts (here the forward) becomes unbounded on $[T-\epsilon,T]$.}
\cite{fleming1985stochastic}.
Since there are several such alternative admissible controls, and each requires infinite energy, comparing these on the basis of the energy that is required is meaningless.
Thus, here, we limit ourselves to a less ambitious goal: we aim to show that an admissible control steering the systems to the desired final marginal exists, %and, moreover, 
that it can be expressed in closed form, and that the resulting controlled evolution coincides with the version of
%shares the same bridge of 
\eqref{system}  conditioned at the initial and final time,
thereby belonging to the same \textit{reciprocal class} as the prior induced by \eqref{system}, cf.\ \cite{Jam1,F2,LK}.

Specifically, we address the following problem.

\begin{probl}\label{Problem2} Provided that the set $\mathcal{U}$ is non-empty, find an explicit construction for an admissible control $u^{*}(\cdot)$ such that i)
$\zeta^{u^{*}}_{t=0}$ is distributed according to $\mu_0$, ii) as $t_n\nearrow T,$ $\zeta^{u^{*}}_{t_n}$ converges in distribution
to a random variable distributed according to $\mu_T$, and iii)
the controlled and uncontrolled evolutions belong to the same reciprocal class\footnote{As noted in the introduction, being in the same reciprocal class means that their probability laws when conditioned on $\zeta_{t=0}=\bar{\zeta}_0, \, \zeta_{t=T}=\bar{\zeta}_T$, for any fixed $\bar{\zeta}_0$, $\bar{\zeta}_T$ in $\R^n$, are identical.}.
\end{probl}

We revisit and discuss issues of optimality in the final section.
Note that, while Assumption \ref{a1} has been removed, Assumption \ref{a2} is always in effect.

\section{Main Results}\label{Results}

Let $\mathcal{N}_0$, $\mathcal{R}_0$ and $\mathcal{N}_T$, $\mathcal{R}_T$ denote the null and range spaces of $\Sigma_0$ and $\Sigma_T$, respectively, and $\Pi_{\mathcal{N}_0}$, $\Pi_{\mathcal{R}_0}$ and 
$\Pi_{\mathcal{N}_T}$, $\Pi_{\mathcal{R}_T}$ denote the corresponding projection operators.
To handle the singularity of the covariances we introduce small perturbations, $\varepsilon_0, \varepsilon_T >0$, and define the perturbed covariances
$$\Sigma_{0,\varepsilon_0}:=\Sigma_0+\varepsilon_0\Pi_{\mathcal{N}_0}, \qquad \Sigma_{T,\varepsilon_T}:=\Sigma_T+\varepsilon_T\Pi_{\mathcal{N}_T}$$ so that $\Sigma_{0,\varepsilon_0},\,\Sigma_{T,\varepsilon_T}\succ 0$ for any $\varepsilon_0, \varepsilon_T>0$. Our goal is to derive the optimal control for our problem as limiting case of the standard results.
Without loss of generality we assume that
$\Sigma_0$ and $\Sigma_T$ are partitioned as
\begin{equation}\label{eq:blocksetting}
\Sigma_0=\begin{bmatrix}
\Lambda_{0} & 0\\
0 & 0
\end{bmatrix}, \qquad \Sigma_T=\begin{bmatrix}
\Lambda_{T} & 0\\
0 & 0
\end{bmatrix}
\end{equation}
and, as a consequence,
$$
\Sigma_{0,\varepsilon_0}=\begin{bmatrix}
\Lambda_{0} & 0\\
0 & \varepsilon_0 I
\end{bmatrix},\qquad \Sigma_{T,\varepsilon_T}=\begin{bmatrix}
\Lambda_{T} & 0\\
0 & \varepsilon_T I
\end{bmatrix}
$$
where $\Lambda_0\in\mathbb{R}^{k_0\times k_0}, \, \Lambda_T\in\mathbb{R}^{k_T\times k_T}$, $\Lambda_0,\Lambda_T\succ 0$ and $0$ and $I$ denote  the zero and the identity matrices of compatible size, respectively.
Indeed, if this were not the case, we can find unitary matrices $U_0,\; U_T$ 
and define a time dependent transformation $U_t$, which can always be chosen to be unitary, as %the group of unitary matrices considered as a topological space is connected{\footnote {\vale Note that the orthogonal group is not connected}}, such that $U(0)=U_0$ and $U(T)=U_T$. 
the set of unitary matrices is path-connected. In fact, one possible choice is 
%$U(t):=\exp\lbrace (1-t)\log U_0+t\log U_T \rbrace$.
$U_t:=\exp[(T-t)/T\log U_0]\exp[t/T\log U_T]$.
Accordingly, changing coordinates in the state space, we return to case \eqref{eq:blocksetting}.

We are now ready to state our first result.\\

\begin{propo}\label{propo_Q}
Let $(P_{\varepsilon_0,\varepsilon_T}(\cdot), \; Q_{\varepsilon_0,\varepsilon_T}(\cdot))$ represent the solution in Proposition \ref{propo_pI} for the system \eqref{system_Lyap}-\eqref{boundary_Lyap} with boundary conditions $\Sigma_{0,\varepsilon_0},\;\Sigma_{T,\varepsilon_T}$. Then
\begin{small}
\begin{align*}
\lim_{\varepsilon_0,\varepsilon_T\rightarrow 0} Q_{\varepsilon_0,\varepsilon_T}&(0)^{-1} = \Phi(T,0)'M(T,0)^{-1}\Phi(T,0)\\
& +(\Sigma_0^{\dagger})^{1/2}(\dfrac{1}{2}I-(\dfrac{1}{4}I+\Sigma_0^{1/2}\hat{\Sigma}_T\Sigma_0^{1/2})^{1/2})(\Sigma_0^{\dagger})^{1/2}\\
& -\Pi_{\mathcal{N}_0}\hat{\Sigma}_T\Sigma_0^{1/2}(\dfrac{1}{2}I+(\dfrac{1}{4}I+\Sigma_0^{1/2}\hat{\Sigma}_T\Sigma_0^{1/2})^{1/2})^{-1}(\Sigma_0^{\dagger})^{1/2}\\
& -(\Sigma_0^{\dagger})^{1/2}(\dfrac{1}{2}I+(\dfrac{1}{4}I+\Sigma_0^{1/2}\hat{\Sigma}_T\Sigma_0^{1/2})^{1/2})^{-1}\Sigma_0^{1/2}\hat{\Sigma}_T\Pi_{\mathcal{N}_0}\\
& -\Pi_{\mathcal{N}_0}\hat{\Sigma}_T\Pi_{\mathcal{N}_0}\\
& +\Pi_{\mathcal{N}_0}\hat{\Sigma}_T\Sigma_0^{1/2}(\dfrac{1}{2}I+(\dfrac{1}{4}I+\Sigma_0^{1/2}\hat{\Sigma}_T\Sigma_0^{1/2})^{1/2})^{-2}\Sigma_0^{1/2}\hat{\Sigma}_T\Pi_{\mathcal{N}_0}
\end{align*}\end{small}
where {\small $\hat{\Sigma}_T:=\Phi(T,0)'M(T,0)^{-1}\Sigma_T M(T,0)^{-1}\Phi(T,0)$}.
Moreover, $$Q(0)^{-1}:=\lim_{\varepsilon_0\rightarrow 0, \, \varepsilon_T\rightarrow 0} Q_{\varepsilon_0,\varepsilon_T}(0)^{-1}$$ has the property that the Riccati differential equation
{\small
\begin{align}\label{Riccati_Q_inv}
\dot{Q}(t)^{-1}=-Q(t)^{-1}A(t)-A(t)'Q(t)^{-1}+Q(t)^{-1}B(t)B(t)'Q(t)^{-1}
\end{align}}
with initial condition $Q(0)^{-1}$ has a finite solution\footnote{ Note that $t=T$ is a {\em finite escape time} for (\ref{Riccati_Q_inv}).}
on $[0,T)$.
\end{propo}
\proof
%%Consider first the limit as $\varepsilon_T\rightarrow0$:
%%\begin{small}
%%\begin{align*}
%%Q_{\varepsilon_0}(0)^{-1}:=\lim_{\varepsilon_T\rightarrow 0} Q_{\varepsilon_0,\varepsilon_T}(0)^{-1} &= \Sigma_{0,\varepsilon_0}^{-1/2}(\dfrac{1}{2}I \\
%%& +\Sigma_{0,\varepsilon_0}^{1/2}\Phi(T,0)'M(T,0)^{-1}\Phi(T,0)\Sigma_{0,\varepsilon_0}^{1/2}\\
%%& -(\dfrac{1}{4}I+\Sigma_{0,\varepsilon_0}^{1/2}\hat{\Sigma}_T\Sigma_{0,\varepsilon_0}^{1/2})^{1/2}) \Sigma_{0,\varepsilon_0}^{-1/2}\\
%%= & \dfrac{1}{2}\Sigma_{0,\varepsilon_0}^{-1}+\Phi(T,0)'M(T,0)^{-1}\Phi(T,0)\\
%%&-\Sigma_{0,\varepsilon_0}^{-1/2}(\dfrac{1}{4}I+\Sigma_{0,\varepsilon_0}^{1/2}\hat{\Sigma}_T\Sigma_{0,\varepsilon_0}^{1/2})^{1/2}\Sigma_{0,\varepsilon_0}^{-1/2} 
%%\end{align*}
%%\end{small}
%where $\hat{\Sigma}_T:=\Phi(T,0)'M(T,0)^{-1}\Sigma_TM(T,0)^{-1}\Phi(T,0)$.\\
To derive the expression for the limit %as $\varepsilon_0\rightarrow0$
we first define the partition
\begin{align*}
\begin{bmatrix}
E_{\varepsilon} & F_{ \varepsilon}\\
F'_{ \varepsilon} & G_{ \varepsilon}
\end{bmatrix}:= \dfrac{1}{2}\Sigma_{0,\varepsilon_0}^{-1}-\Sigma_{0,\varepsilon_0}^{-1/2}(\dfrac{1}{4}I+\Sigma_{0,\varepsilon_0}^{1/2}\hat{\Sigma}_{T, \varepsilon_T}\Sigma_{0,\varepsilon_0}^{1/2})^{1/2}\Sigma_{0,\varepsilon_0}^{-1/2}
\end{align*}
where $\hat{\Sigma}_{T,\varepsilon_T}:=\Phi(T,0)'M(T,0)^{-1}\Sigma_{T,\varepsilon_T} M(T,0)^{-1}\Phi(T,0)$, $E_{\varepsilon}\in\mathbb{R}^{k_0\times k_0}$, $F_{\varepsilon}\in \mathbb{R}^{k_0\times(n-k_0)}$ and $G_{\varepsilon}\in\mathbb{R}^{(n-k_0)\times (n-k_0)}$. By rearranging terms we get:
\begin{align*}
\Bigg(\Sigma_{0,\varepsilon_0}^{1/2}\begin{bmatrix}
E_{\varepsilon} & F_{\varepsilon}\\
F'_{\varepsilon} & G_{\varepsilon}
\end{bmatrix} \Sigma_{0,\varepsilon_0}^{1/2} -\dfrac{1}{2}I \Bigg)^{2}&=\dfrac{1}{4}I+\Sigma_{0,\varepsilon_0}^{1/2}\hat{\Sigma}_{T, \varepsilon_T}\Sigma_{0,\varepsilon_0}^{1/2}
\end{align*}
and by developing the square and simplifying:
\begin{align}\label{befor_lim}
\begin{bmatrix}
E_{\varepsilon} & F_{ \varepsilon}\\
F'_{ \varepsilon} & G_{ \varepsilon}
\end{bmatrix} \Sigma_{0,\varepsilon_0}\begin{bmatrix}
E_{\varepsilon} & F_{\varepsilon}\\
F'_{\varepsilon} & G_{\varepsilon}
\end{bmatrix}-\begin{bmatrix}
E_{\varepsilon} & F_{\varepsilon}\\
F'_{\varepsilon} & G_{\varepsilon}
\end{bmatrix}=\hat{\Sigma}_{T, \varepsilon_T}.
\end{align}
Let $\hat{\Sigma}_T:=\lim_{\varepsilon_T \rightarrow 0 }\,\hat{\Sigma}_{T,\varepsilon_T}$, $E:=\lim_{ \varepsilon_0 \rightarrow0,\varepsilon_T \rightarrow 0 }\, E_{ \varepsilon}$, $F:=\lim_{ \varepsilon_0 \rightarrow0,\varepsilon_T \rightarrow 0}\, F_{ \varepsilon}$ and $G:=\lim_{ \varepsilon_0 \rightarrow0,\varepsilon_T \rightarrow 0}\, G_{\varepsilon}$.
Then, taking the limit for $\varepsilon_0\rightarrow0$, $\varepsilon_T\rightarrow0$ of \eqref{befor_lim} and suitably partitioning $\hat{\Sigma}_T$ we obtain the following system of equations in $E$, $F$ and $G$: 
\begin{align*}
\begin{bmatrix}
E & F\\
F' & G
\end{bmatrix}  
\begin{bmatrix}
\Lambda_{0} & 0 \\
0 & 0
\end{bmatrix}
\begin{bmatrix}
E & F\\
F' & G
\end{bmatrix}-\begin{bmatrix}
E & F\\
F' & G
\end{bmatrix}=\begin{bmatrix}\hat{\Sigma}_T^{E} & \hat{\Sigma}_T^{F} \\ \hat{\Sigma}_T^{F'} & \hat{\Sigma}_T^{G}\end{bmatrix}
\end{align*}
or equivalently
\begin{align*}
\begin{cases}
E\Lambda_{0}E -E=\hat{\Sigma}_T^{E}\\
E\Lambda_{0}F -F=\hat{\Sigma}_T^{F}\\
F'\Lambda_{0}F-G=\hat{\Sigma}_T^{G}
\end{cases}.
\end{align*}
Only the first equation is quadratic and has as solutions:
$$E_\pm =\Lambda_{0}^{-1/2}\Big(\pm \Big(\dfrac{1}{4}I+\Lambda_{0}^{1/2}\hat{\Sigma}_T^{E}\Lambda_{0}^{1/2}\Big)^{1/2}+\dfrac{1}{2}I \Big)\Lambda_{0}^{-1/2}.$$
To resolve the choice of the sign we observe that $E_{\varepsilon}$ is continuous as $\varepsilon_0 \rightarrow 0,\, \varepsilon_T \rightarrow 0$ thus, by comparison, the minus sign is the consistent one.
Therefore, the solutions for $F$ and $G$ are:
\begin{align*}
F&=-\Lambda_{0}^{-1/2}\Big(\Big(\dfrac{1}{4}I+\Lambda_{0}^{1/2}\hat{\Sigma}_T^{E}\Lambda_{0}^{1/2}\Big)^{1/2}+\dfrac{1}{2}I \Big)^{-1}\Lambda_{0}^{1/2}\hat{\Sigma}_T^{F}\\
G &= -\hat{\Sigma}_T^{G} 
+ \hat{\Sigma}_T^{F'}\Lambda_{0}^{1/2}\Big( \Big(\dfrac{1}{4}I+\Lambda_{0}^{1/2}\hat{\Sigma}_T^{E}\Lambda_{0}^{1/2}\Big)^{1/2}+\dfrac{1}{2}I \Big)^{-2}\Lambda_{0}^{1/2}\hat{\Sigma}_T^{F}.
\end{align*}
Then the
expression for $\lim_{\varepsilon_0\rightarrow 0, \varepsilon_T\rightarrow 0}Q_{\varepsilon_0,\varepsilon_T}(0)^{-1}$ in the statement can be verified utilizing the definitions of $\Sigma_0$ and $\Pi_{\mathcal{N}_0}$.

Next, to show that the Riccati differential equation \eqref{Riccati_Q_inv} has finite well defined solution in $[0,T)$ we introduce the following change of coordinates:
$$x(t)=N(T,0)^{-1/2}\Phi(0,t)\zeta(t)$$
where $N(t_1,t_0):=\int_{t_0}^{t_1}\Phi(t_0,\tau)B(\tau)B(\tau)'\Phi(t_0,\tau)'d\tau$ is the controllability Gramian which is non singular for all $t_0<t_1$, $t_0,t_1\, \in [0,T]$.
In this new coordinates the dynamics simplify to
$$dx(t)=\underbrace{N(T,0)^{-1/2}\Phi(0,t)B(t)}_{=:B_{\text{new}}}dW_t$$
for the newly defined diffusion matrix $B_{\text{new}}$. This ensures that for the new reachability Gramian it holds that $M_{\text{new}}(T,0)=I$.\\
Accordingly under the new coordinates \eqref{system_Lyap_P} becomes $\dot{Q}_{\text{new}}(t)=-B_{\text{new}}(t)B_{\text{new}}(t)'$ and the relation between $Q_{\text{new}}(t)$ and $Q(t)$ is given by
$$Q_{\text{new}}(t):=N(T,0)^{-1/2}\Phi(0,t)Q(t)\Phi(0,t)'N(T,0)^{-1/2}.$$
It follows that
\begin{align*}
Q_{\text{new}}(t)^{-1}=&(Q_{\text{new}}(0)-M_{\text{new}}(t,0))^{-1}\\
=&-M_{\text{new}}(t,0)^{-1}-M_{\text{new}}(t,0)^{-1}(Q_{\text{new}}(0)^{-1} \\ 
& -M_{\text{new}}(t,0)^{-1})^{-1}M_{\text{new}}(t,0)^{-1}.
\end{align*}
Therefore, showing that the Riccati differential equation \eqref{Riccati_Q_inv} admits finite solution in $[0,T)$ reduces to showing that 
\begin{align}\label{non_sing}
\det\left(Q_{\text{new}}(0)^{-1}-M_{\text{new}}(t,0)^{-1}\right)\neq 0 \quad \forall t\in[0,T).
\end{align}
The expression $(Q_{\text{new}}(0)^{-1}-M_{\text{new}}(t,0)^{-1})$ is  maximal for $t=T$; 
here ``maximal'' is to be interpreted in the sense of the natural partial order on positive semi-definite matrices. 
Moreover, the following relation holds:
\[
M(T,0)^{-1}= \Phi(0,T)'N(T,0)^{-1}\Phi(0,T).
\]
Therefore, for all $t<T$, we have that
\begin{align}\label{upper_bound}
Q_{\text{new}}(0)^{-1}& -M_{\text{new}}(t,0)^{-1} \notag \\
& = I + N(T,0)^{\frac{1}{2}}\begin{bmatrix}
E & F \\ F' & G
\end{bmatrix} N(T,0)^{\frac{1}{2}} - M_{\text{new}}(t,0)^{-1} \notag \\
& \prec  N(T,0)^{\frac{1}{2}}\begin{bmatrix}
E & F \\ F' & G
\end{bmatrix} N(T,0)^{\frac{1}{2}}.
\end{align}
We will now establish that the following holds
\begin{align}\label{neg_sd}
\begin{bmatrix}
E & F \\ F' & G
\end{bmatrix} \preceq 0.
\end{align}
which in view of \eqref{upper_bound} implies \eqref{non_sing}.

To this end, we make the following observations:
\begin{itemize}
\item[a)]the matrix $E$ is negative definite as $$\frac{1}{2}I-\left(\frac{1}{4}I+\Lambda_{0}^{1/2}\hat{\Sigma}_T^{E}\Lambda_{0}^{1/2}\right)^{1/2} \prec 0;$$
\item[b)] by definition $\hat{\Sigma}_T\succeq 0$ therefore $\hat{\Sigma}_T^{G}\succeq 0$ and by taking the (generalized) Schur complement $$\hat{\Sigma}_T^{G}-\hat{\Sigma}_T^{F'}(\hat{\Sigma}_T^{E})^{\dagger}\,\hat{\Sigma}_T^{F}\succeq 0.$$
\end{itemize}
From a) $(-E)\succ 0$, therefore proving \eqref{neg_sd} is equivalent to proving that
\begin{align*}
(-G)-(-F')(-E)^{-1}(-F)\succeq 0.
\end{align*}
As a consequence of b), this amounts to show that
\begin{align}\label{ineq}
F'\Lambda_{0}F-F'E^{-1}F \preceq \hat{\Sigma}_T^{F'}(\hat{\Sigma}_T^{E})^{\dagger}\hat{\Sigma}_T^{F}.
\end{align}
We define $\Theta:=\left(\frac{1}{4}I+\Lambda_{0}^{1/2}\hat{\Sigma}_T^{E}\Lambda_{0}^{1/2}\right)^{1/2}$ and we observe that $\Lambda_{0}^{-1/2}(\hat{\Sigma}_T^{E})^{\dagger}\Lambda_{0}^{-1/2}=\left(\Theta^{2}-\frac{1}{4}I \right)^{-1}$.
Then \eqref{ineq} reduces to
%\begin{align*}
%-F'\Sigma_{0,+}F+F'E^{-1}F= -\hat{\Sigma}_T^{F'}\Sigma_{0,+}^{1/2}\Big\lbrace \left(\Theta+\frac{1}{2}I\right)^{-2}-\left(\Theta+\frac{1}{2}I\right)^{-1}\left(\frac{1}{2}I-\Theta\right)^{-1}\left(\Theta+\frac{1}{2}I\right)^{-1} \Big\rbrace\Sigma_{0,+}^{1/2}\hat{\Sigma}_{T}^{F}
%\end{align*}
\begin{footnotesize}
\[
\left(\Theta+\frac{1}{2}I\right)^{-2}-\left(\Theta+\frac{1}{2}I\right)^{-1}\left(\frac{1}{2}I-\Theta\right)^{-1}\left(\Theta+\frac{1}{2}I\right)^{-1} \preceq \left( \Theta^{2}-\frac{1}{4}I \right)^{-1}\]
%\[\left(\Theta+\frac{1}{2}I\right)^{-1}\left[ I - \left(\frac{1}{2}I-\Theta\right)^{-1} \right]\left(\Theta+\frac{1}{2}I \right)^{-1} \? \left( \Theta^{2}-\frac{1}{4}I \right)^{-1}\]
%\[\left[ I - \left(\frac{1}{2}I-\Theta\right)^{-1} \right]\? \left(\Theta+\frac{1}{2}I\right)\left( \Theta^{2}-\frac{1}{4}I \right)^{-1}\left(\Theta+\frac{1}{2}I\right)\]
%\[\left[ I - \left(\frac{1}{2}I-\Theta\right)^{-1} \right]\? \left(\Theta+\frac{1}{2}I\right)\left( \Theta-\frac{1}{2}I \right)^{-1}\left(\Theta+\frac{1}{2}I\right)^{-1}\left(\Theta+\frac{1}{2}I\right)\]
%\[I\?\left(\Theta+\frac{1}{2}I-I \right)\left(\Theta-\frac{1}{2}I\right)^{-1}\]
\end{footnotesize}\black
which further reduces to
\begin{footnotesize}
\[\left(\Theta-\frac{1}{2}I\right) \preceq \left(\Theta-\frac{1}{2}I\right).
\]
\end{footnotesize}\black
Therefore \eqref{ineq} holds with equality, completing the proof. 
\qed

\begin{rem}
By Proposition \ref{propo_Q}, $Q(t)$ is non-singular in $[0,T)$.
\end{rem}

Next, we state the analogous result for $P(\cdot)$.

\begin{propo}\label{propo_P}
Let $(P_{\varepsilon_0,\varepsilon_T}(\cdot), \; Q_{\varepsilon_0,\varepsilon_T}(\cdot))$ be as defined in Proposition \eqref{propo_Q}. Then, in view of Proposition \ref{rem_on_P(T)}
\begin{small}
\begin{align*}
\lim_{\varepsilon_0\rightarrow 0, \, \varepsilon_T\rightarrow 0}& P_{\varepsilon_0,\varepsilon_T}(T)^{-1}=  M(T,0)^{-1}\\
& +(\Sigma_T^{\dagger})^{1/2}(\dfrac{1}{2}I -(\dfrac{1}{4}I+\Sigma_T^{1/2}\hat{\Sigma}_0\Sigma_T^{1/2})^{1/2})(\Sigma_T^{\dagger})^{1/2}\\
& - \Pi_{\mathcal{N}_T}\hat{\Sigma}_0\Sigma_T^{1/2}(\dfrac{1}{2}I+(\dfrac{1}{4}I+\Sigma_T^{1/2}\hat{\Sigma}_0\Sigma_T^{1/2})^{1/2})^{-1}(\Sigma_T^{\dagger})^{1/2}\\
& -(\Sigma_T^{\dagger})^{1/2}(\dfrac{1}{2}I+(\dfrac{1}{4}I+\Sigma_T^{1/2}\hat{\Sigma}_0\Sigma_T^{1/2})^{1/2})^{-1}\Sigma_T^{1/2}\hat{\Sigma}_0\Pi_{\mathcal{N}_T}\\
& - \Pi_{\mathcal{N}_T}\hat{\Sigma}_0\Pi_{\mathcal{N}_T} \\ 
& + \Pi_{\mathcal{N}_T}\hat{\Sigma}_0\Sigma_T^{1/2}(\dfrac{1}{2}I+(\dfrac{1}{4}I+\Sigma_T^{1/2}\hat{\Sigma}_0\Sigma_T^{1/2})^{1/2})^{-2}\Sigma_T^{1/2}\hat{\Sigma}_0\Pi_{\mathcal{N}_T}
\end{align*}\end{small}
where $\hat{\Sigma}_0:=M(T,0)^{-1}\Phi(T,0)\Sigma_0 \Phi(T,0)'M(T,0)^{-1}$.
Moreover, the limit $$P(T)^{-1}:=\lim_{\varepsilon_0\rightarrow 0, \, \varepsilon_T\rightarrow 0} P_{\varepsilon_0,\varepsilon_T}(T)^{-1}$$ has the property that the Riccati differential equation
\begin{align*}
\dot{P}(t)^{-1}=-P(t)^{-1}A(t)-A(t)'P(t)^{-1}-P(t)^{-1}B(t)B(t)'P(t)^{-1}
\end{align*}
with final condition $P(T)^{-1}$ has finite %well-defined 
solution in $(0,T]$. 
\end{propo}
The proof of Proposition \ref{propo_P} is completely symmetric to the proof of Proposition \ref{propo_Q}.

\begin{rem}
By Proposition \ref{propo_P} $P(t)$ is non-singular in $(0,T]$.
\end{rem}

\begin{propo} \label{block_struct}
Under the considered coordinate basis, $P(0)$ and $Q(T)$ have the following block diagonal structure
\begin{align}\label{eq_block_struct}
P(0)=\Pi_{\mathcal{R}_0}P(0)\Pi_{\mathcal{R}_{0}}=\begin{bmatrix}
P_{0,+} & 0 \\ 0 & 0
\end{bmatrix},
\\
Q(T)=\Pi_{\mathcal{R}_T}Q(T)\Pi_{\mathcal{R}_T}=\begin{bmatrix}
Q_{T,+} & 0 \\
0 & 0
\end{bmatrix}
\end{align}
where $P_{0,+}\in\mathbb{R}^{k_0 \times k_0}$ and $Q_{T,+}\in\mathbb{R}^{k_T \times k_T}$ are both non-singular.
\end{propo}
The proof of Proposition \ref{block_struct} is deferred to the Appendix.

\begin{thm}\label{main_thm} Define $\hat{u}(t):=-B(t)'Q(t)^{-1}\zeta_t$ and consider the controlled evolution \eqref{controlled_ev} for the choice $u(t)=\hat{u}(t)$. Then, the solution to the linear stochastic differential equation
\begin{equation}\label{controled_main_thm}
d\zeta^{\hat{u}}_t=\hat{A}(t)\zeta^{\hat{u}}_tdt+B(t)dW_t,\quad 0\leq t < T
\end{equation}
with $\zeta^{\hat{u}}_{t=0}=\zeta_0$ a.s. and $\hat{A}(t):= A(t)-B(t)B(t)'Q(t)^{-1}$ %converges in $L^{2}(\Omega,\mathcal{F},\mathbb{P})$
is such that its mean function is zero and, as $t\nearrow T$, its variance  converges to $\Sigma_T$.
\end{thm}
\proof
It can be checked that $\mathbb{E}\lbrace \zeta^{\hat{u}}_t\rbrace=0$, since $\E\{\zeta_0\}=0$.
To verify the second-order statistics, we define $\Sigma(t):=\mathbb{E}\lbrace \zeta^{\hat{u}}_t{\zeta_t^{\hat{u}}}^\prime\rbrace$. Then, $\Sigma(t)$ satisfies the Lyapunov differential equation
\begin{align}\label{Lyap_cov}
\dot{\Sigma}(t)=\hat{A}(t)\Sigma(t)+\Sigma(t)\hat{A}(t)'+B(t)B(t)' 
\end{align}
with initial condition $\Sigma(0)=\mathbb{E}\lbrace \zeta_0\zeta_0' \rbrace=\Sigma_0 .$
It can be easily verified that the candidate solution:
\begin{align}\label{candidate_sol}\Sigma(t)=(P(t)^{-1}+Q(t)^{-1})^{-1}
\end{align} 
satisfies the differential equation \eqref{Lyap_cov} in the open interval $(0,T)$, where $Q(\cdot)^{-1}$ and $P(\cdot)^{-1}$ are as defined in Propositions \ref{propo_Q} and \ref{propo_P}.
We claim, in addition, that relation \eqref{candidate_sol}  holds at the end points $t=0$ and $t=T$ in the sense that 
\begin{subequations}\label{boundary}
\begin{align}\label{t_0}
P(0)(I+Q(0)^{-1}P(0))^{-1}=\Sigma_0\\
\label{t_T}
Q(T)(P(T)^{-1}Q(T)+I)^{-1}=\Sigma_T.
\end{align}
\end{subequations}
%since $P(0)$ and $Q(T)$ are singular.\\
Let us focus on \eqref{t_0}. %We want to show that $\Sigma_{\varepsilon_0}(t)$ is continuous at $\varepsilon_0=0,\, t=0$.\\
By simple algebra, $(Q(t)^{-1}+P(t)^{-1})^{-1}=P(t)(I+Q(t)^{-1}P(t))^{-1}$ for $t\neq 0$. Then, the limit %for $t\searrow 0$ the following holds
$$\lim_{t\searrow 0}P(t)(I+Q(t)^{-1}P(t))^{-1}=P(0)(I+Q(0)^{-1}P(0))^{-1}$$
is well-defined provided that $\det(I+Q(0)^{-1}P(0))\neq 0$.

To show that $(I+Q(0)^{-1}P(0))$ is indeed not singular we recall that, in view of Proposition \ref{block_struct},  $P(0)$ has the block-diagonal structure in \eqref{eq_block_struct}.
Then, conformally partitioning
$$Q(0)^{-1}=\begin{bmatrix}
{E} & {F}\\
{F}' & {G}
\end{bmatrix},$$
the quantity $(I+Q(0)^{-1}P(0))$ can be expressed as
$$\begin{bmatrix}
I &0\\
0 & I
\end{bmatrix} +\begin{bmatrix}
{E} & {F}\\
{F}' & {G}
\end{bmatrix}\begin{bmatrix}
P_{0,+} & 0\\
0 & 0
\end{bmatrix}= 
\begin{bmatrix}
I+{E}P_{0,+} & 0\\ 
{F}'P_{0,+} & I
\end{bmatrix},
$$
whose determinant is $\det(I+{E}P_{0,+})$. We want to show now that this quantity is non zero.

We recall that, for $\varepsilon_0>0$, the results in Theorem \ref{Thm_part_I} apply and the following condition holds
\begin{align}\label{epsilon}
P_{\varepsilon_0}(0)(I+Q_{\varepsilon_0}(0)^{-1}P_{\varepsilon_0}(0))^{-1}=\Sigma_{0,\varepsilon_0}.\end{align} 
We observe that $Q_{\varepsilon_0}(0),\, Q_{\varepsilon_0}(0)^{-1}$ and $P_{\varepsilon_0}(0)$ are well defined as $\varepsilon_0\rightarrow 0 $ %$, \, t=0$
 and we introduce the representation $Q_{\varepsilon_0}^{-1}(0)=Q(0)^{-1}+\mathcalO(\varepsilon_0)$, $P_{\varepsilon_0}(0)=P(0)+\mathcalO(\varepsilon_0)$ so \eqref{epsilon} gives
 \begin{footnotesize}
\begin{align*}
&
\bigg(\begin{bmatrix}
P_{0,+} & 0\\
0 & 0
\end{bmatrix} \bhs+\mathcalO(\varepsilon_0) \bigg)\left(I +\bigg(\begin{bmatrix}
{E} & {F}\\
{F}' & {G}
\end{bmatrix}\bhs+\mathcalO(\varepsilon_0) \bigg)\bigg(\begin{bmatrix}
P_{0,+} & 0\\
0 & 0
\end{bmatrix} \bhs+\mathcalO(\varepsilon_0) \bigg) \right)^{-1}
\\ 
&\hspace*{5cm}=\begin{bmatrix}
\Sigma_{0,+} & 0\\
0 & \varepsilon_0
\end{bmatrix},
\end{align*}
\end{footnotesize}
%$$\begin{bmatrix}
%\Sigma_{0,+} & 0\\
%0 & \varepsilon_0
%\end{bmatrix}\left(\bigg(\begin{bmatrix}
%P_{0,+} & 0\\
%0 & 0
%\end{bmatrix} +\mathcalO(\varepsilon_0)\bigg)^{-1}+\bigg(\begin{bmatrix}
%\tilde{E} & \tilde{F}\\
%\tilde{F}' & \tilde{G}
%\end{bmatrix}+\mathcalO(\varepsilon_0)\bigg)\right)=I$$
%$$\begin{bmatrix}
%\Sigma_{0,+} & 0\\
%0 & \varepsilon_0
%\end{bmatrix}\left(I +\bigg(\begin{bmatrix}
%\tilde{E} & \tilde{F}\\
%\tilde{F}' & \tilde{G}
%\end{bmatrix}+\mathcalO(\varepsilon_0)\bigg)\bigg(\begin{bmatrix}
%P_{0,+} & 0\\
%0 & 0
%\end{bmatrix} +\mathcalO(\varepsilon_0)\bigg)\right)= \bigg(\begin{bmatrix}
%P_{0,+} & 0\\
%0 & 0
%\end{bmatrix} +\mathcalO(\varepsilon_0)\bigg) $$
which further reduces to
\begin{small}
$$\begin{bmatrix}
\Sigma_{0,+} & 0\\
0 & \varepsilon_0
\end{bmatrix}\left( 
\begin{bmatrix}
I+{E}P_{0,+} & 0\\
{F}'P_{0,+} & I 
\end{bmatrix} +\mathcalO(\varepsilon_0)
\right)= \begin{bmatrix}
P_{0,+} & 0\\
0 & 0
\end{bmatrix} +\mathcalO(\varepsilon_0).
$$
\end{small}
The limit as $\varepsilon_0\searrow 0$
is well defined and given by
\begin{small}
$$\begin{bmatrix}
\Sigma_{0,+}(I+{E}P_{0,+}) & 0\\
0 & 0
\end{bmatrix} =\begin{bmatrix}
P_{0,+} & 0\\
0 & 0
\end{bmatrix}$$
\end{small}from which it follows that $\det(I+{E}P_{0,+})\neq 0$.

Moreover, we observe that
\begin{align*}
\lim_{t\searrow 0}(Q(t)^{-1}+P(t)^{-1})^{-1}&=\lim_{t\searrow 0}P(t)(I+Q(t)^{-1}P(t))^{-1} \\
&=P(0)(I+Q(0)^{-1}P(0))^{-1}\\
&=\lim_{\varepsilon_0\searrow 0}P_{\varepsilon_0}(0)(I+Q_{\varepsilon_0}(0)^{-1}P_{\varepsilon_0}(0))^{-1}\\
&=\lim_{\varepsilon_0\searrow 0}\Sigma_{0,\varepsilon_0}=\Sigma_0
\end{align*}
and \eqref{t_0} follows.
Equation \eqref{t_T} can be shown by using a symmetric argument.
%\begin{align*}
%\lim_{t\nearrow T}(Q(t)^{-1}+P(t)^{-1})^{-1}&=\lim_{t\nearrow T}Q(t)(P(t)^{-1}Q(t)+I)^{-1}\\
%&= Q(T)(P(T)^{-1}Q(T)+I)^{-1}\\
%&= \lim_{\varepsilon_T\searrow 0}Q_{\varepsilon_T}(T)(P_{\varepsilon_T}(T)^{-1}Q_{\varepsilon_T}(T)+I)^{-1}\\
%&= \lim_{\varepsilon_T\searrow 0} \Sigma_{T,\varepsilon_T}= \Sigma_T.
%\end{align*}
%%Finally, since $\lim_{t\nearrow T}(Q(t)^{-1}+P(t)^{-1})^{-1}=\Sigma_T$, we conclude that $\zeta^{\hat{u}}_t$ converges in $L^{2}(\Omega,\mathcal{F},\mathbb{P})$ to a random variables distributed according to $\mu_T$ which concludes the proof.\\
\qed

Consider now Problem \ref{Problem2} in the case when $A(t)\equiv 0$, $B(t)\equiv I_n$ and the two marginals are delta masses at $0$. Then, we readily get $Q(t)=(1-t)I_n$ and  $\hat{u}(t):=-\frac{1}{1-t}\zeta_t$. Thus, the Brownian Bridge of Section \ref{BrownianB} is indeed the solution of (the generalized Schr\"odinger Bridge)  Problem \ref{Problem2}. 

%Finally, we make the following remark.
%\begin{rem}The uncontrolled evolution \eqref{system} and the controlled evolution \eqref{controlled_ev} share the same bridge, that is, the probability law of \eqref{system} and when conditioned on $\zeta(0)=\bar{\zeta}_0$, $\zeta(T)=\bar{\zeta}_T$, for any fixed $\bar{\zeta}_0$, $\bar{\zeta}_T$ \textcolor[rgb]{0,.4,.5}{in $\R^n$}, is the same as the probability law of \eqref{controlled_ev} when conditioned on the same event.\end{rem}We refer to the analogous result in \cite{chen2016optimal} for a proof of this fact.

%\section{Time symmetry of the solution}\label{time-symmetry}

The following result basically amounts to the fact that the uncontrolled evolution (prior) and the bridge evolution belong to the same reciprocal class.
\begin{propo}
The uncontrolled evolution \eqref{system} and the controlled evolution \eqref{controled_main_thm} share the same bridge, that is, the probability law of \eqref{system}  when conditioned on $\zeta_{t=0}=\bar{\zeta}_0$, $\zeta_{t=T}=\bar{\zeta}_T$, for any fixed $\bar{\zeta}_0$, $\bar{\zeta}_T$ in $\R^n$, is the same as the probability law of \eqref{controled_main_thm} when conditioned on the same event.
\end{propo}

\proof
Let $R_1$ and $R_2$ be, respectively, the solutions of
$$\dot{R}_1(t)=A(t)R_1(t)+R_1(t)A(t)'-B(t)B(t)'$$
$$\dot{R}_2(t)=\hat{A}(t)R_2(t)+R_2(t)\hat{A}(t)'-B(t)B(t)'$$
with terminal conditions $R_1(T)=R_2(T)=0$. Then
a proof of the statement follows by the same argument of \cite[Theorem 11]{chen2016optimal}
observing that $\lim_{t\rightarrow T}R_1(t)Q(t)^{-1}R_2(t)=0$.
\qed

We have seen
above that the forward drift  field $b_+(\zeta,t)$ of the bridge is obtained from the uncontrolled one by an additive perturbation of the form $-B(t)B(t)'Q(t)^{-1}\zeta$ which makes the forward drift  $b_+(\zeta^{\hat{u}}(t),t)$ unbounded as $t\nearrow T$.  The bridge process solving Problem \ref{Problem2} also admits a reverse time differential whose drift field $b_-(\zeta,t)$ is $b_+(\zeta,t)+B(t)B(t)'\Sigma(t)^{-1}\zeta$ \cite[Lemma 2.3]{BLP1979}. Here $\Sigma(t)$ is the covariance of the optimal $\zeta^{\hat{u}}(t)$. Because of (\ref{candidate_sol}), we now get
\begin{align}\label{backward}
b_-(\zeta_t,t)&=b_+(\zeta_t,t)+B(t)B(t)'\Sigma(t)^{-1}\zeta_t \notag \\
&=\left[A(t)-B(t)B(t)'Q(t)^{-1}+B(t)B(t)'\Sigma(t)^{-1}\right]\zeta_t \notag \\
&=\left[A(t)+B(t)B(t)'P(t)^{-1}\right]\zeta_t.
\end{align}
From (\ref{backward}), we see that, in a specular way, $b_-(\zeta^{\hat{u}}_t,t)$ is unbounded as $t\searrow 0$. For instance, in the case of the Brownian bridge considered at the end of the previous section, we have $P(t)=tI_n$ and the backward drift is $b_-(\zeta^{\hat{u}}_t,t)=\frac{1}{t}\zeta^{\hat{u}}_t$. 

\section{Numerical Example{\ignore s}} \label{Numerical_ex}
\begin{figure}[htb]
\begin{minipage}[b]{0.47\textwidth}
	\centering
\includegraphics[clip,width=1\textwidth]{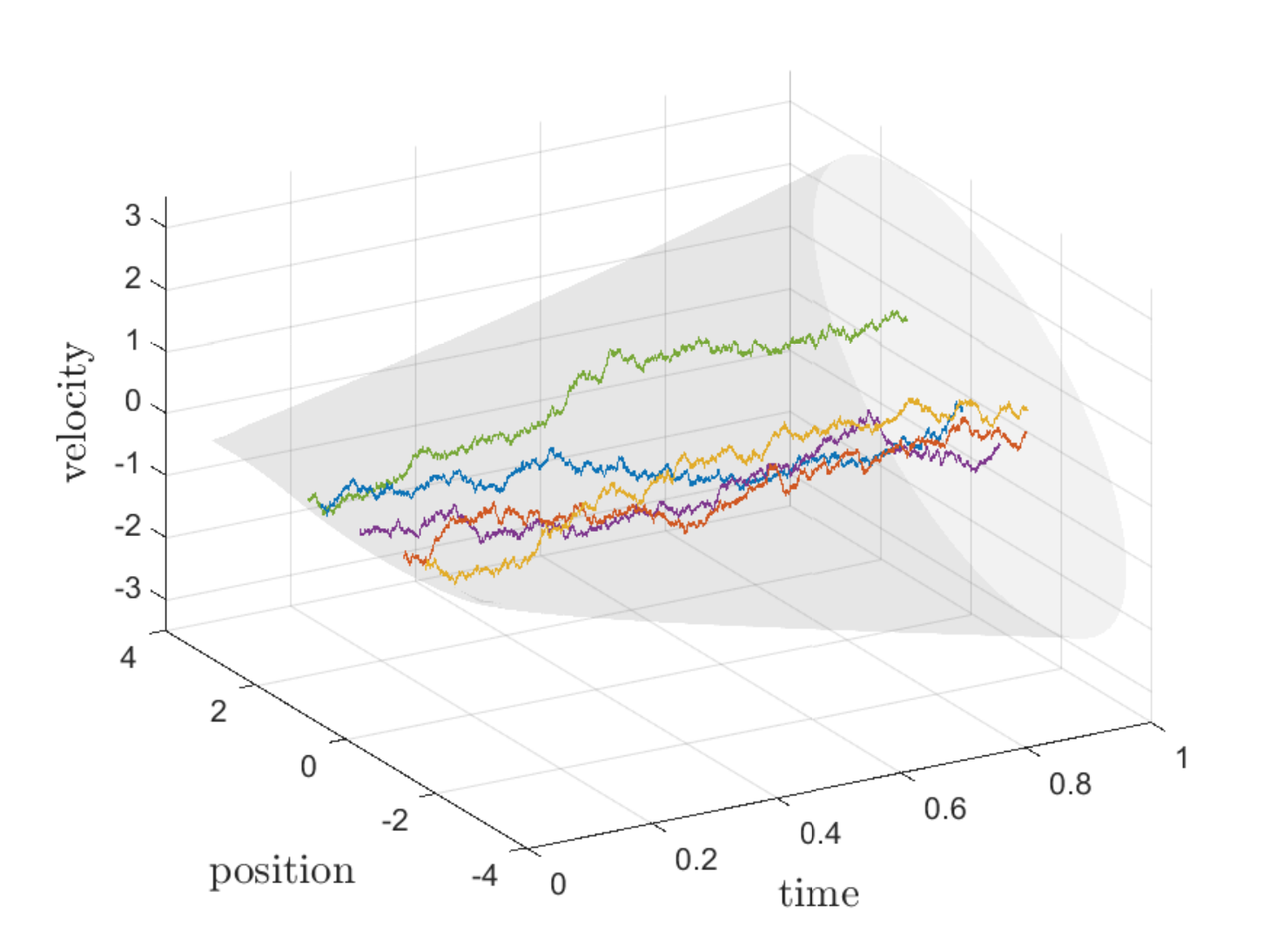}
	\caption{\small \ignore{Example 1:} Sample paths for the uncontrolled evolution.}
\label{fig1}
\end{minipage}
\begin{minipage}[b]{0.47\textwidth}
%\begin{figure}
	\centering
\includegraphics[clip,width=1\textwidth]{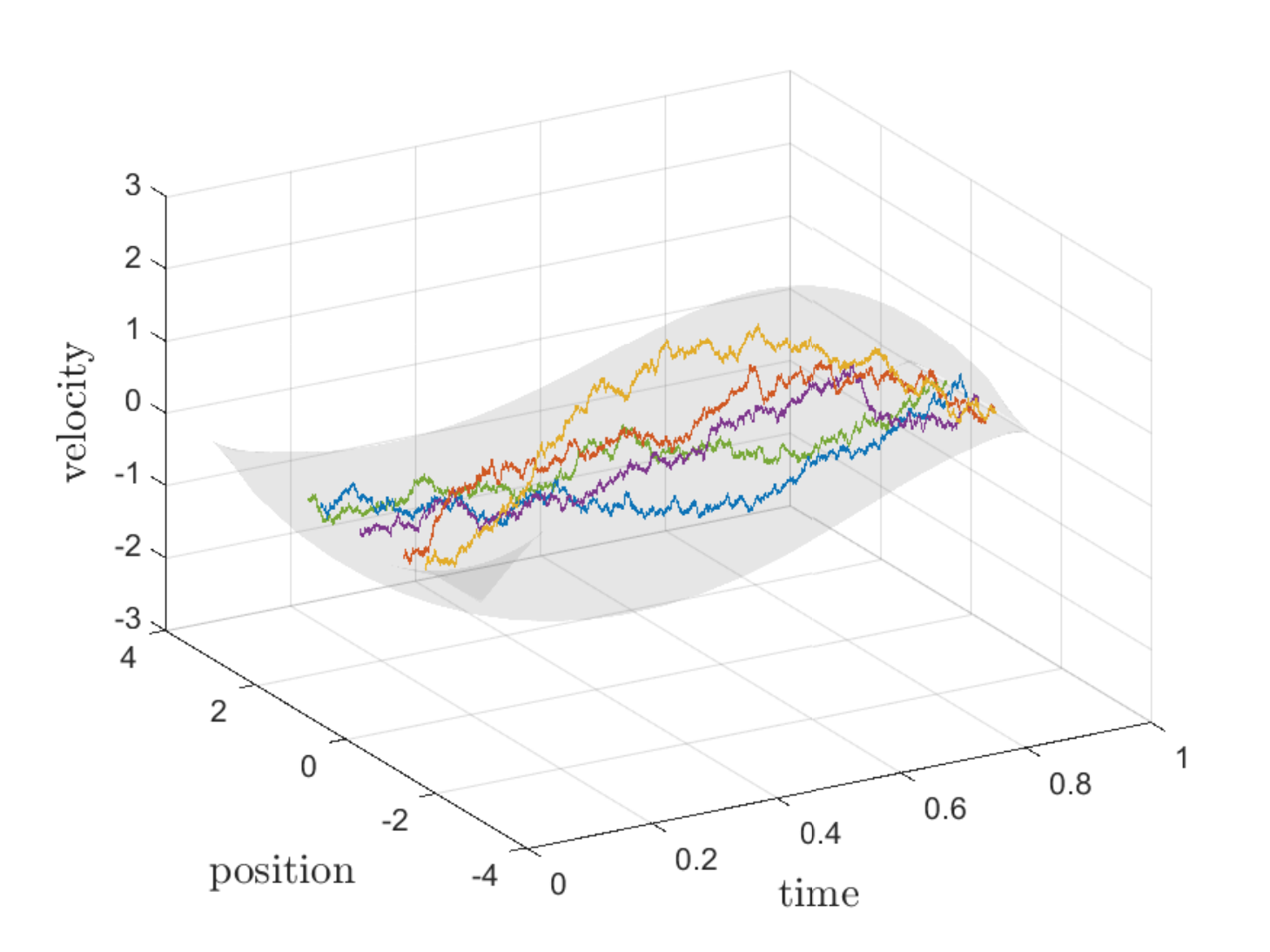}
	\caption{\small \ignore{Example 1:} Sample paths for the controlled evolution.}
\label{fig2}
\end{minipage}
\end{figure}

\begin{figure}[htb]
\begin{minipage}[b]{0.47\textwidth}
	\centering
\includegraphics[clip,width=1\textwidth]{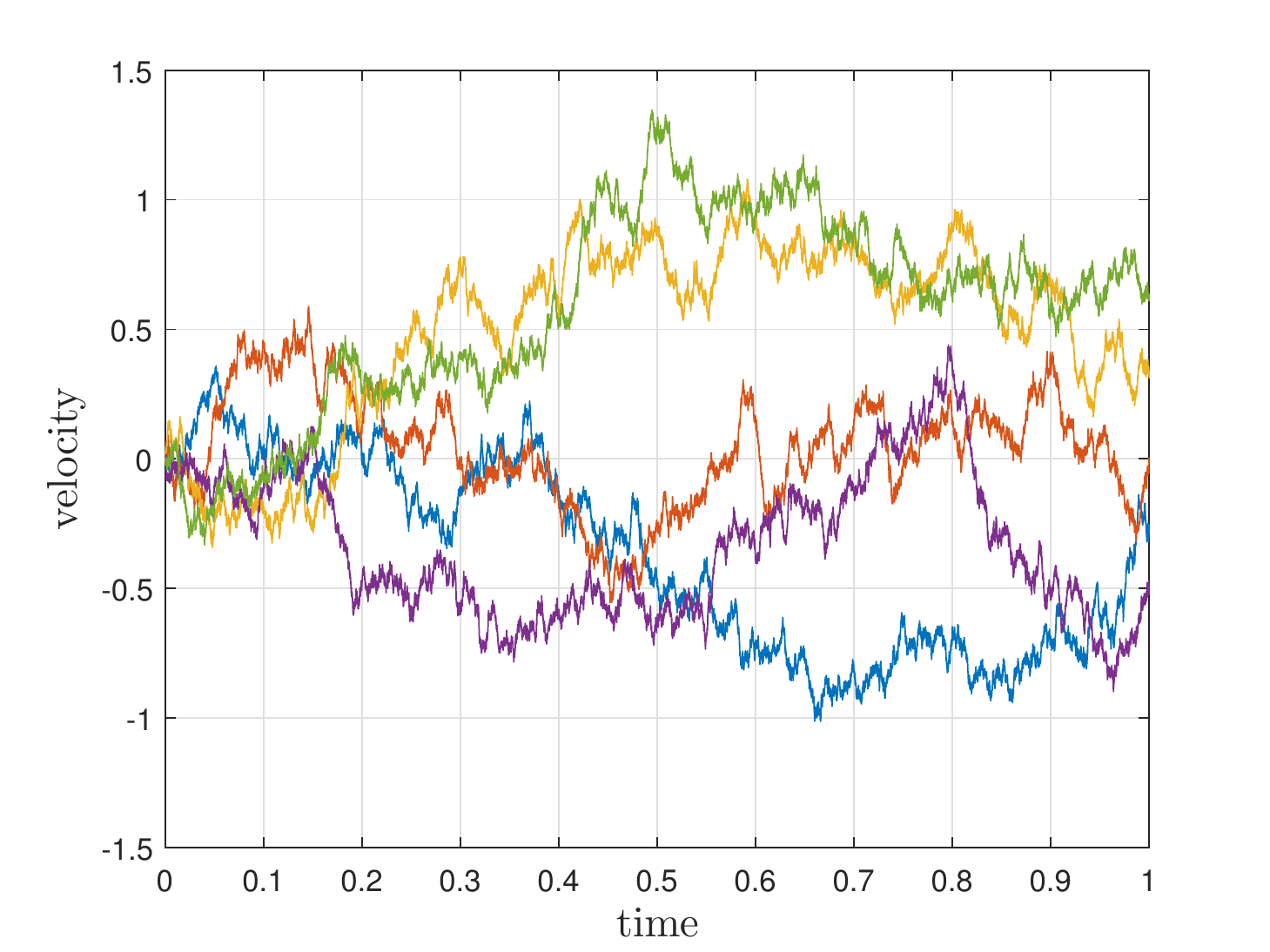}
	\caption{\small \ignore{Example 1:} Sample paths for the uncontrolled velocity.}
\label{fig3}
\end{minipage}
\begin{minipage}[b]{0.47\textwidth}
	\centering
\includegraphics[clip,width=1\textwidth]{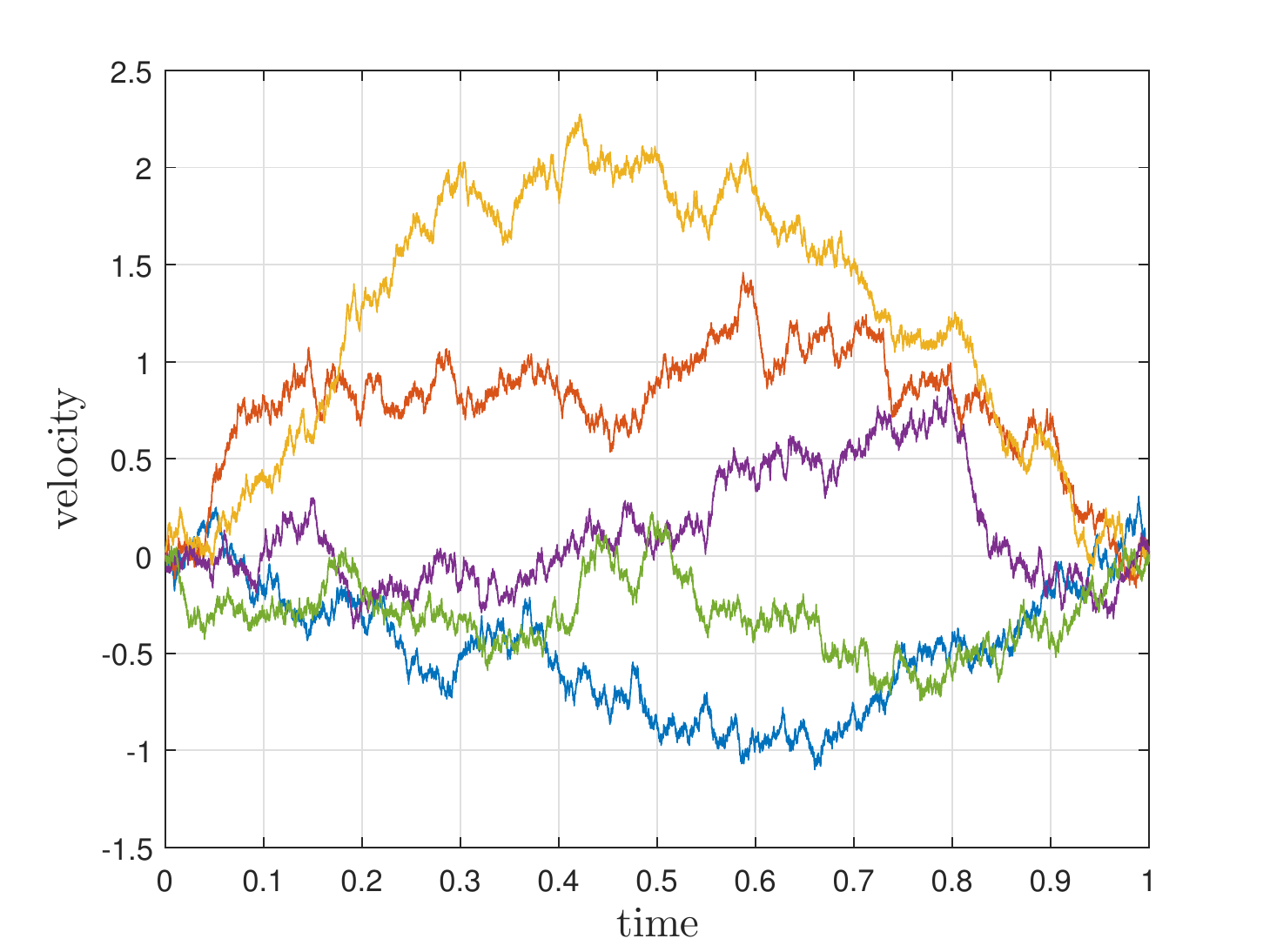}
	\caption{\small \ignore{Example 1:} Sample paths for the controlled velocity.}
\label{fig4}
\end{minipage}
\end{figure}

\begin{figure}[htb]
\begin{minipage}[b]{0.47\textwidth}
	\centering
\includegraphics[clip,width=1\textwidth]{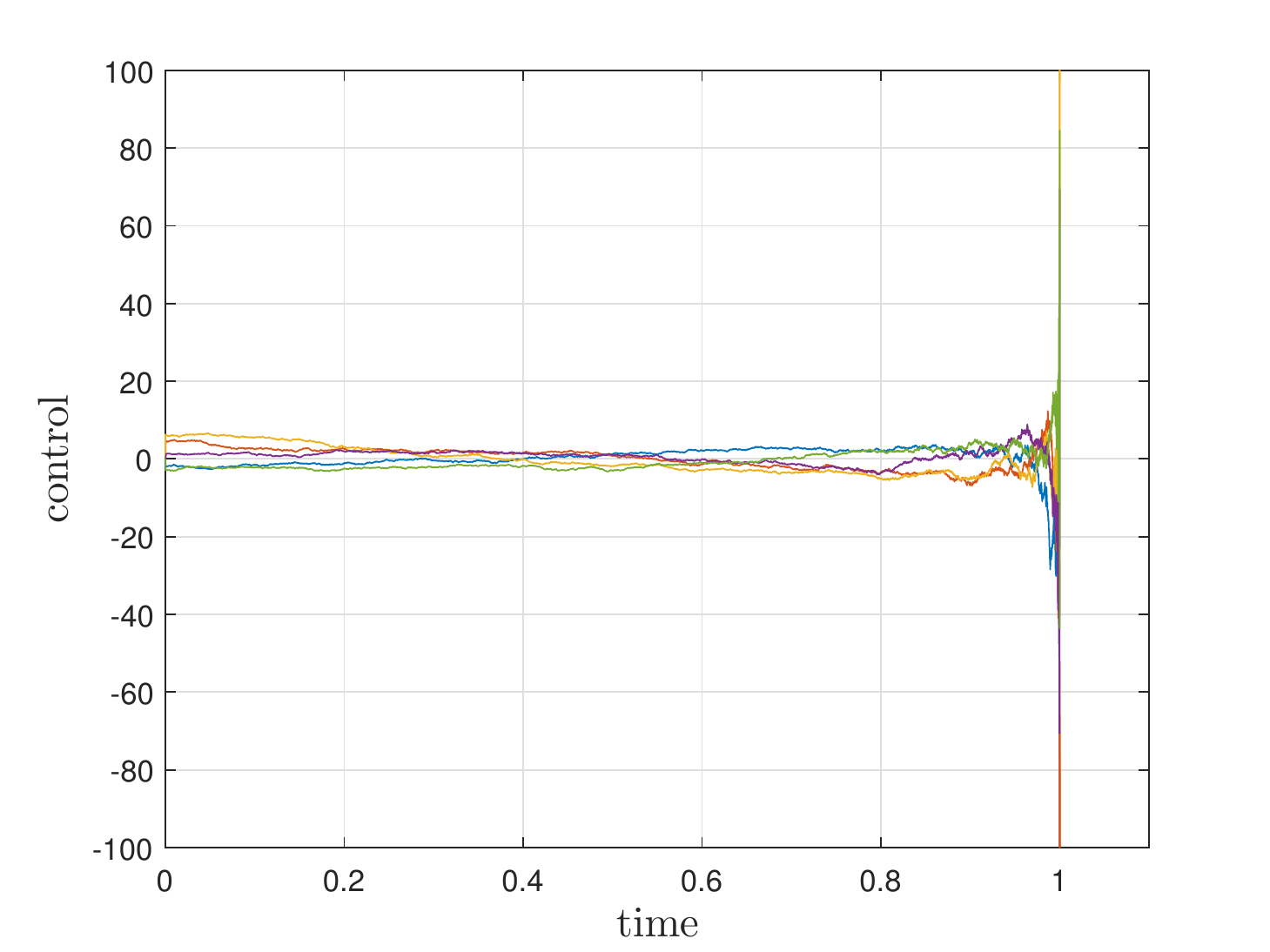}
	\caption{\small \ignore{Example 1:} Control inputs.}
\label{fig4.51}
\end{minipage}
\ignore{
\begin{minipage}[b]{0.47\textwidth}
	\centering
\includegraphics[clip,width=1\textwidth]{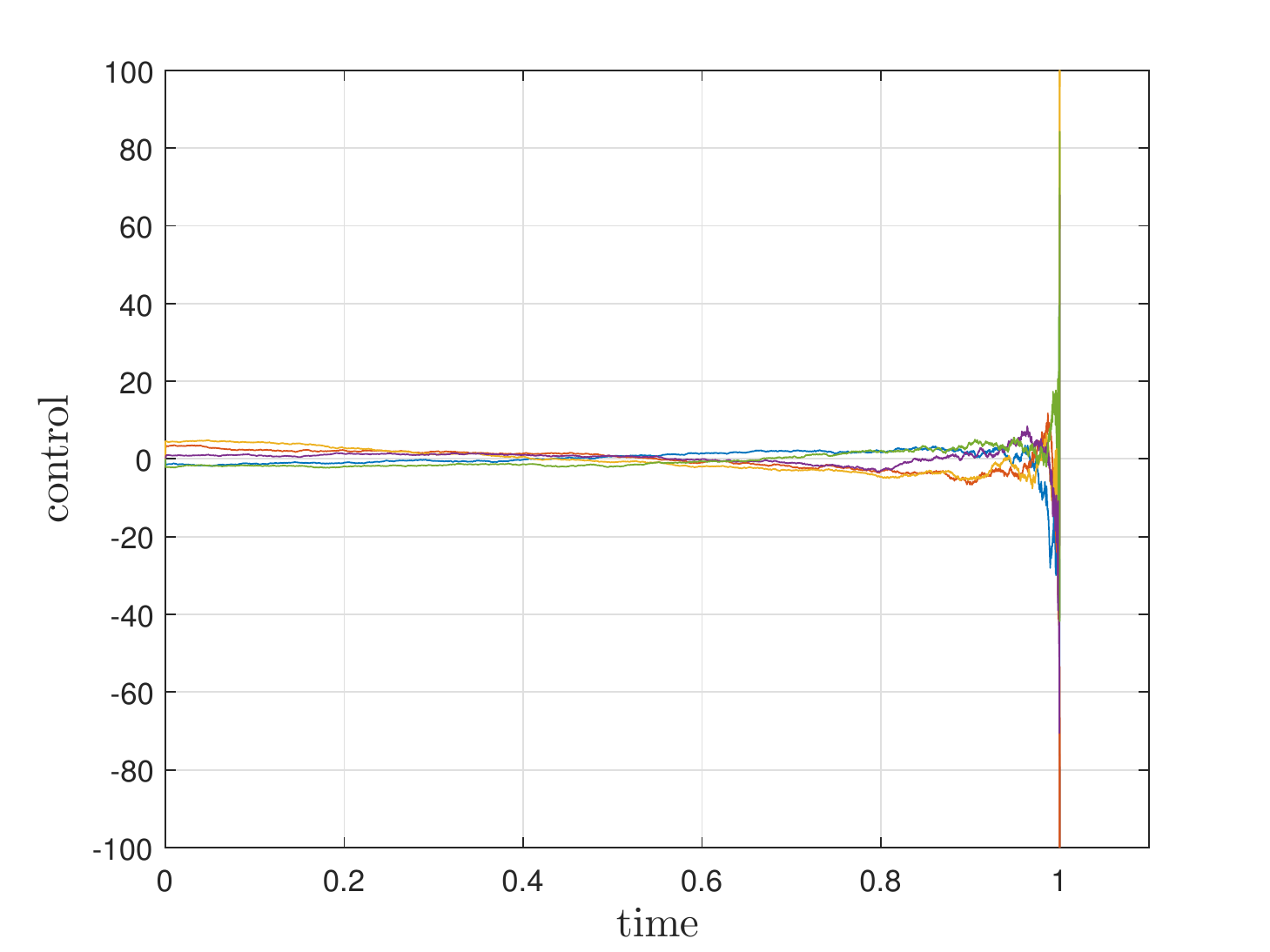}
	\caption{\small Example 2: Control inputs.}
\label{fig4.52}
\end{minipage}
}
\end{figure}

\ignore{\textbf{Example 1.}} To illustrate the results, we consider inertial particles subject to random acceleration:
%As an academic example for illustrative purposes, we consider the following stochastic velocity model:
$$\begin{bmatrix}
dx(t) \\ dv(t)
\end{bmatrix}=\begin{bmatrix}
0 & 1 \\ 0 & 0
\end{bmatrix} \begin{bmatrix}
x(t) \\ v(t)
\end{bmatrix}dt+
\begin{bmatrix}
0 \\ 1
\end{bmatrix} u(t)dt + \begin{bmatrix}
0 \\ 1
\end{bmatrix}dW_t
$$
where $x(t)$ represents the position, $v(t)$ the velocity and $u(t)$ the control force. We address the problem of  steering a cloud of particles, obeying such a dynamics, from an initial distribution with $x(0)\sim \mathcal{N}(0,1), v(0)\sim \delta_0$, to a final one with $x(T)\sim\mathcal{N}(0,0.2),v(T)\sim \delta_0$, where $T=1$. Position and velocity are independent at the initial and final time. Figures \ref{fig1} and \ref{fig2} display the sample paths in the phase space $(x,v)$ as a function of time for the uncontrolled and controlled evolution, respectively, whereas figures \ref{fig3} and \ref{fig4} display the uncontrolled velocity and the controlled velocity, respectively. The shaded regions in the phase plots represent the ``$3\sigma$'' confidence region. The optimal feedback control $u(t)$ has been derived according to the results in Proposition \ref{propo_Q}.
Figure \ref{fig4.51} displays the corresponding control inputs diverging, as expected, at  $t=1$.

\ignore{

\textbf{Example 2.} Next, we consider a modification of the above model where we account for the presence of a viscosity term and of an harmonic potential in the velocity equation:
$$\begin{bmatrix}
dx(t) \\ dv(t)
\end{bmatrix}=\begin{bmatrix}
0 & 1 \\ -1 & -1
\end{bmatrix} \begin{bmatrix}
x(t) \\ v(t)
\end{bmatrix}dt+
\begin{bmatrix}
0 \\ 1
\end{bmatrix} u(t)dt + \begin{bmatrix}
0 \\ 1
\end{bmatrix}dW_t.
$$
Observe that now the velocity evolution is no more Markovian. As before, we address the problem of  steering a cloud of particles obeying such dynamics, from an initial distribution with $x(0)\sim \mathcal{N}(0,1)$, $v(0)\sim \delta_0$,  to a final one with $x(1)\sim\mathcal{N}(0,0.2)$, $v(1)\sim \delta_0$. As before position and velocity are assumed to be independent at the initial and final time. Figure \ref{fig4.52} displays the control inputs derived according to the results in Proposition \ref{propo_Q}. Figures \ref{fig5} and \ref{fig6} display the sample paths in the phase space $(x,v)$ as a function of time for the uncontrolled and controlled evolution, respectively, whereas figures \ref{fig7} and \ref{fig8} display the uncontrolled velocity and the Brownian bridge for controlled velocity, respectively. As above, the shaded regions in the phase plots represent the ``$3\sigma$'' confidence region.

\begin{figure}[htb]
\begin{minipage}[b]{0.47\textwidth}
	\centering
\includegraphics[clip,width=1\textwidth]{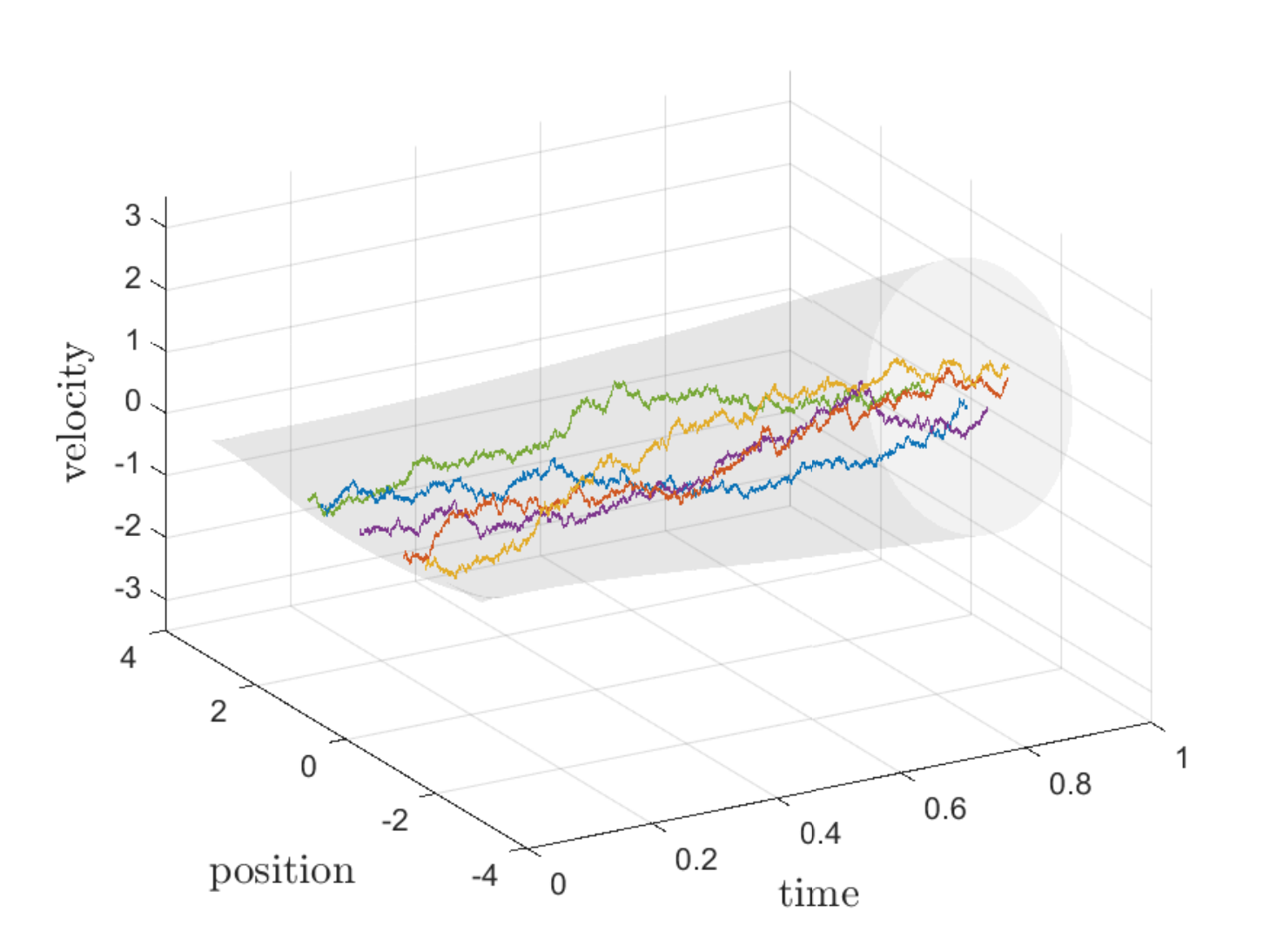}
	\caption{\small Example 2: Sample paths for the uncontrolled evolution.}
\label{fig5}
\end{minipage}
\begin{minipage}[b]{0.47\textwidth}
%\begin{figure}
	\centering
\includegraphics[clip,width=1\textwidth]{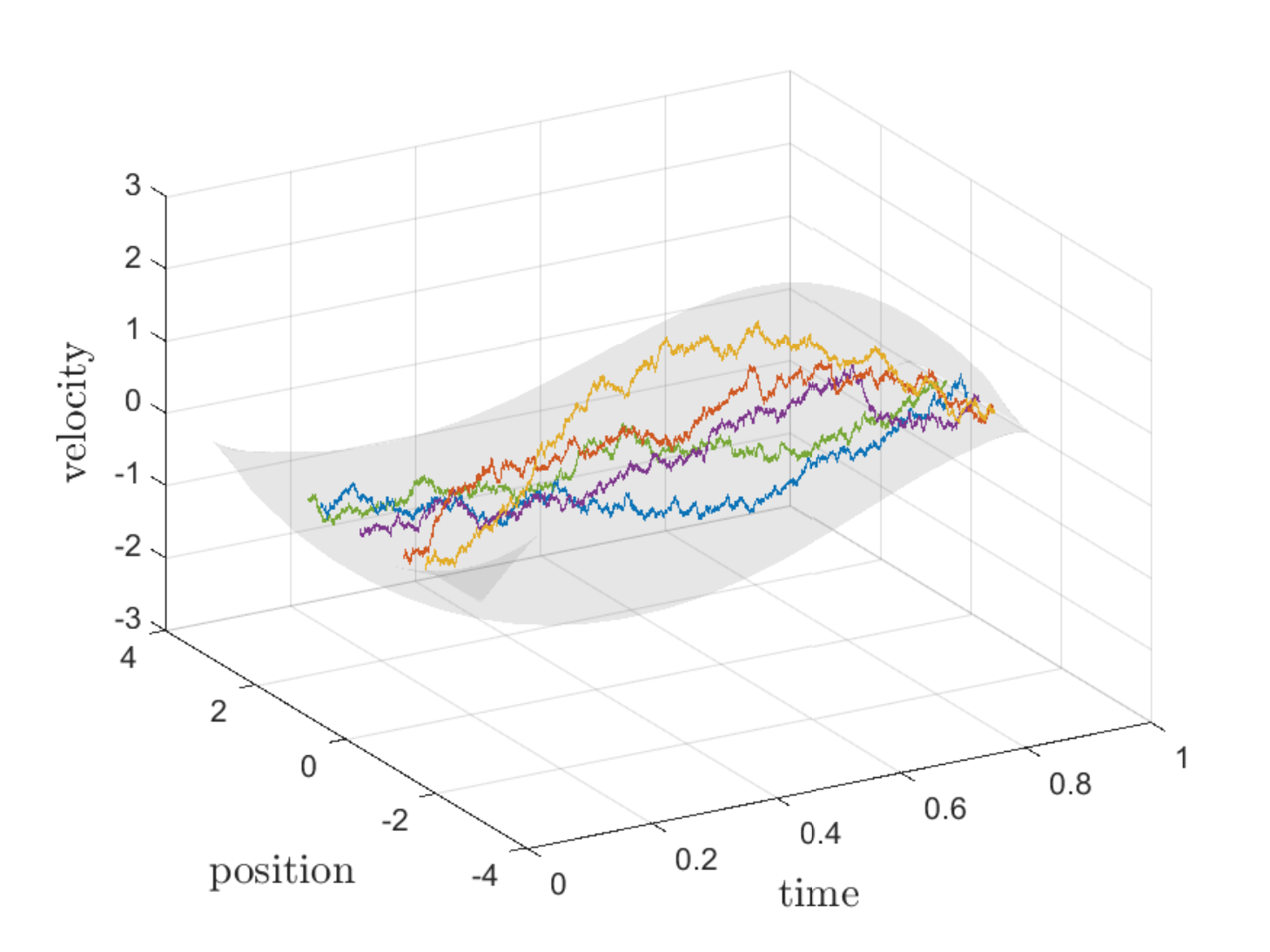}
	\caption{\small Example 2: Sample paths for the controlled evolution.}
\label{fig6}
\end{minipage}
\end{figure}

\begin{figure}[htb]
\begin{minipage}[b]{0.47\textwidth}
	\centering
\includegraphics[clip,width=1\textwidth]{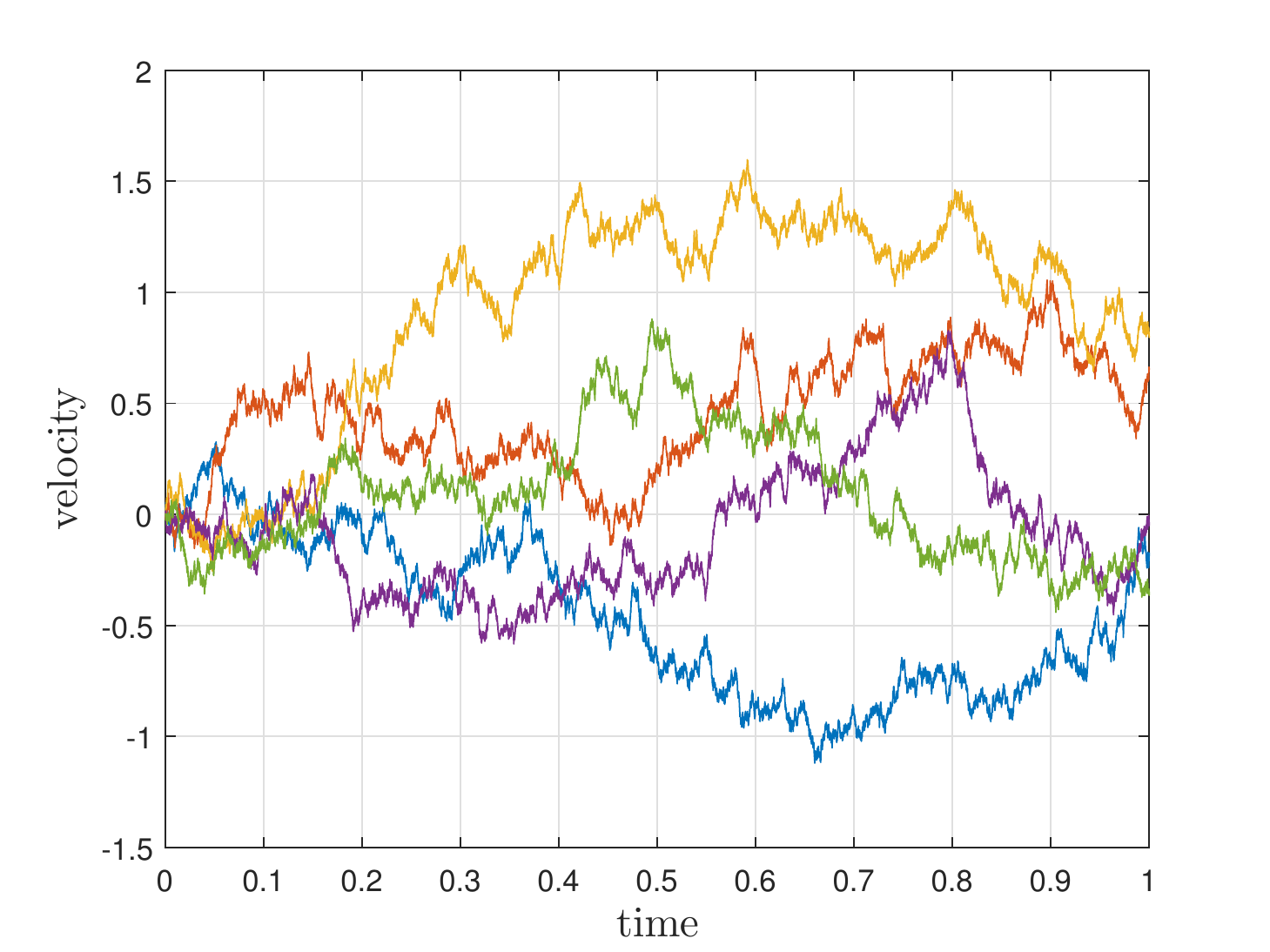}
	\caption{\small Example 2: Sample paths for the uncontrolled velocity.}
\label{fig7}
\end{minipage}
\begin{minipage}[b]{0.47\textwidth}
	\centering
\includegraphics[clip,width=1\textwidth]{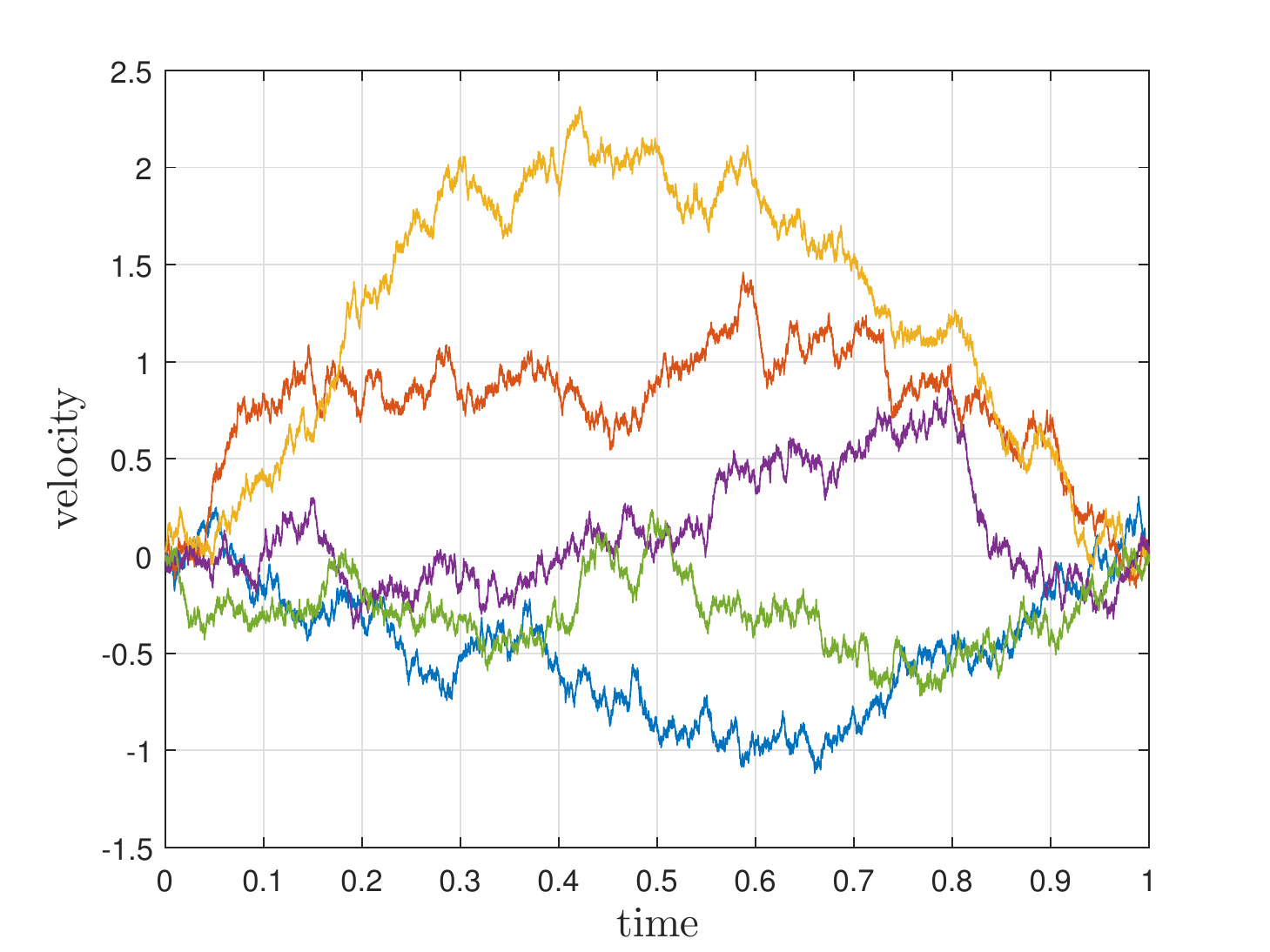}
	\caption{\small Example 2: Sample paths for the controlled velocity.}
\label{fig8}
\end{minipage}
\end{figure}
}

\section{Concluding remarks and future directions}\label{Conclusion}
We have addressed the problem of steering a linear stochastic system between two end-point
{\em degenerate} 
Gaussian distributions. 
Specifically, 
we have shown that there exists an admissible control steering the system to the desired final configuration, and that this can be expressed in closed form. Moreover, the resulting controlled evolution
belongs to the same reciprocal class of the uncontrolled one.

The issue of ``optimality'' is rather subtle as, necessarily, any control strategy steering the system to the desired degenerate marginal grows unbounded at $t=T$. This same issue arose in the case of fixed boundary conditions in \cite{fleming1985stochastic}, i.e., Dirac marginals,  and was addressed by  controlling $\zeta_t$ over the smaller time interval $T-\varepsilon$ while introducing a suitable penalty term $F(\zeta_{T-\varepsilon})$.

An alternative strategy, in the spirit of what has already been presented in the present paper, is to define ``optimality'' of the control indirectly by expressing it as the limit of a standard (non-degenerate) case for a suitable vanishing perturbation (brought in to remove the degeneracy) and requiring that the limit is independent of the chosen perturbation. 
%$\varepsilon_0\Pi_0$, $\varepsilon_T\Pi_T$ such that $\Sigma_0+\varepsilon_o\Pi_0, \, \Sigma_T+\varepsilon_t\Pi_T\succ 0$ and denote by $u^{\ast}_\varepsilon$ the corresponding minimum energy control obtained as in 
%Theorem \ref{Thm_part_I}. Then, we may define $u^{\ast}$ to be optimal for our problem if $u^{\ast}=\lim_{\varepsilon\rightarrow 0}u^{\ast}_\varepsilon$ and the limit is independent of the chosen perturbation. 

Interestingly, in the special case where degeneracy appears at one of the terminal points, {\em a natural choice for the optimal control is inherited by the time-symmetry of the problem}. Specifically, when the starting marginal at $t=0$ is degenerate, and therefore $\Sigma_0$ singular while $\Sigma_T$ is not, one would not perceive any issue when comparing controls. Indeed, in the forward time-direction admissible controls only need to have finite energy and thereby can be compared directly. When the degeneracy is reversed, with $\Sigma_T$ being singular and $\Sigma_0$ not, then the forward-in-time controls are unbounded but the backwards-in-time are not and can be compared for deciding on the optimal one and the corresponding probability law on paths. The probability law on the paths does not give a preference to a time-arrow and can be viewed either way. Thus, it induces a unique choice for the corresponding infinite-energy control in the forward-in-time direction. Optimality of the forward-in-time control can then be claimed on this basis (optimality of the law).

Comparing control strategies in general, when the energy required is infinite, appears to require a deeper exploration and it is hoped that it will be revisited in future work.

%In this case the optimal control still grows unbounded at the final time $T$, but the problem may be solved backward (for non singular initial covariance) by standard approaches.
%Finally, a third possibility concerns the restricted case in which only the terminal covariance is singular. In this case the optimal control still grows unbounded at the final time $T$, but the problem may be solved backward (for non singular initial covariance) by standard approaches.

%A {\color{blue} deeper study of these strategies will be the } subject of future investigation.

\begin{appendix}
\subsection{Proof of Proposition \ref{rem_on_P(T)}}

\proof As in \cite{chen2016optimal} we introduce the following change of variables:
$$x(t)=N(T,0)^{-1/2}\Phi(0,t)\xi_t$$
%$$dx(t)= \underbrace{N(T,0)^{-1/2}\Phi(0,t)B(t)}_{=:B_{new}}dW_t$$
where $N(\cdot,\cdot)$ is the controllability Gramian of the system.
Under this new coordinates system we have that 
$$
\tilde{\Sigma}_{0,\text{new}}=N(T,0)^{-1/2}\tilde{\Sigma}_0 N(T,0)^{-1/2},
$$
$$
\tilde{\Sigma}_{T,\text{new}}=N(T,0)^{-1/2}\Phi(0,T)\tilde{\Sigma}_{T}\Phi(0,T)'N(T,0)^{-1/2},
$$
$$P_{\text{new}}(t)=N(T,0)^{-1/2}\Phi(0,t)P(t)\Phi(0,t)'N(T,0)^{-1/2}.$$
Then, by the very same argument in \cite{chen2016optimal}, we obtain the following two sets of final conditions for the system  
of Lyapunov differential equations \eqref{system_Lyap}:
\begin{small}
\begin{align*}
P^{\pm}(T)=& \Phi(T,0)N(T,0)^{1/2} 
\tilde{\Sigma}_{T,\text{new}}^{1/2}
\bigg(\tilde{\Sigma}_T+\frac{1}{2}I\pm \\
& \bhs\bhs\bhs \left(\tilde{\Sigma}_{T,\text{new}}^{1/2}\tilde{\Sigma}_{0,\text{new}}\tilde{\Sigma}_{T,\text{new}}^{1/2}+\frac{1}{4}I\right)^{1/2}
\bigg)^{-1}
\tilde{\Sigma}_{T,\text{new}}^{1/2}
N(T,0)^{1/2}\Phi(T,0)'\\
Q^{\pm}(T)&=(\tilde{\Sigma}_T^{-1}-P^{\pm}(T)^{-1})^{-1}
\end{align*}
\end{small}
and again by mimicking \cite{chen2016optimal} it can be shown that $P^{-}(t)$ and $Q^{-}(t)$ are non-singular on $[0,T]$.

To recover the formula in the statement let us define $G:=\tilde{\Sigma}_T^{1/2}\Phi(0,T)'N(T,0)^{-1/2}$ so that $\tilde{\Sigma}_{T,\text{new}}=G'G$. Then, there exists an orthogonal matrix $U$ such that $UG=G'U'$, $UG\succ 0$. Thus, $UG$ is the unique symmetric positive definite square root of $\tilde{\Sigma}_{T,\text{new}}$. Substituting %in \eqref{PT_A} 
in the expression for $P^{-}(T)$ and rearranging terms,
\begin{small}
\begin{align*} P(T)= & \Big( \Phi(0,T)'N^{-1}\Phi(0,T)
\\ & +\frac{1}{2}\Sigma_T^{-1}-\Phi(0,T)'N^{-1/2}(UG)^{-1}\Big(
UGN^{-1/2}\Sigma_0N^{-1/2}G'U' \\
& +\frac{1}{4}I \Big)(G'U')^{-1}N^{-1/2}\Phi(0,T) \Big), \end{align*}
\end{small}\black
where $N:=N(T,0)$. 
Finally, to recover the formula in the statement it is enough to observe that $M(T,0)=\Phi(T,0)N(T,0)\Phi(T,0)'$ and that the following equality holds:
\begin{small}
\begin{align*} 
U'\Big(& UGN^{-1/2}\Sigma_0N^{-1/2}G'U'+\frac{1}{4}I \Big)^{1/2}U \\
& =\Big(\Sigma_T^{1/2}M(T,0)^{-1}\Phi(T,0)\Sigma_0\Phi(T,0)'M(T,0)^{-1}\Sigma_T^{1/2}+\frac{1}{4}I\Big)^{1/2}.
\end{align*}
\end{small}
\qed
\subsection{Proof of Proposition \ref{block_struct}}
\proof We show the result for $P(0)$; the result for $Q(T)$ is completely symmetric.

For $\varepsilon_0>0$ Theorem \ref{Thm_part_I} applies and the following holds
\begin{align}\label{from_PartI}
\Sigma_{0,\varepsilon_0}(I+Q_{\varepsilon_0}(0)^{-1}P_{\varepsilon_0}(0))=P_{\varepsilon_0}(0).
\end{align}
Moreover, we recall that in view of Proposition \ref{propo_Q} $\lim_{\varepsilon_0\rightarrow 0}Q_{\varepsilon_0}(0)^{-1}$ is well defined. We partition $P_{\varepsilon_0}(0)$ and $Q_{\varepsilon_0}(0)^{-1}$ conformally to $\Sigma_0$ as
$$ P_{\varepsilon_0}(0):= 
\begin{bmatrix}
 R+\mathcalO(\varepsilon_0) & S+\mathcalO(\varepsilon_0)\\
 S'+\mathcalO(\varepsilon_0) & V+\mathcalO(\varepsilon_0)
\end{bmatrix},
\; 
Q_{\varepsilon_0}(0)^{-1}:=
\begin{bmatrix}
E_{\varepsilon_0} & F_{\varepsilon_0}\\
F_{\varepsilon_0}' & G_{\varepsilon_0}
\end{bmatrix}.
$$
and by substituting in \eqref{from_PartI}  and taking the limit we obtain
$$\begin{bmatrix}
\Lambda_{0}& 0 \\ 0 & 0
\end{bmatrix}\left( 
\begin{bmatrix}
I & 0 \\ 0 & I
\end{bmatrix}
+
\begin{bmatrix}
E & F \\ F' & G
\end{bmatrix}
\begin{bmatrix}
R & S \\ S' & V
\end{bmatrix}
\right)=
\begin{bmatrix}
R & S \\ S' & V
\end{bmatrix}
$$
where $\lim_{\varepsilon_0\rightarrow 0} E_{\varepsilon_0}= E$, $\lim_{\varepsilon_0\rightarrow 0} F_{\varepsilon_0}= F$ and $\lim_{\varepsilon_0\rightarrow 0} G_{\varepsilon_0}= G$.\\ 
Then, by solving the associated system of equations in $R, \, S$ and $V$,  we get that $S=V=0$ and that
\begin{align}\label{E_pd}
\Lambda_{0}+\Lambda_{0}ER=R.
\end{align}
Finally, from \eqref{E_pd} it is immediate to see that $R$ is non-singular as $\Lambda_{0}=(I-\Lambda_0E)R$ is positive definite by definition.
\qed
\end{appendix}

%\subsection*{Acknowledgement}
%Partial support was provided by NSF under grants %1665031,
%1807664, 1839441 (TTG) and 1901599 (YC), and by AFOSR under grant FA9550-17-1-0435 (TTG), and by the University of Padova Research Project CPDA 140897 (MP).\\

%%\spacingset{0.95}

%\bibliographystyle{plain} abbrv
\bibliographystyle{abbrv}
\bibliography{biblio}

\end{document}